\documentclass[preprint,5p,twocolumn]{elsarticle}
\graphicspath{ {./images/} }

\usepackage{lipsum}
\usepackage{forest}
\usepackage{amsfonts}
\usepackage{graphicx}
\usepackage{epstopdf}
\usepackage{algorithmic}

\usepackage{commath,amsopn}
\usepackage{multirow}
\usepackage{hyperref}
\usepackage{stmaryrd}
\usepackage{xifthen}
\usepackage{booktabs}
\usepackage{graphicx}
\usepackage{subfigure}
\usepackage{epstopdf}
\usepackage{chemformula,cleveref}
\usepackage[switch,columnwise]{lineno}

\definecolor{myred}{rgb}{0.7,0.1,0.45}
\definecolor{myrvs}{rgb}{0.25,0.45,0.85}
\definecolor{whucolor}{rgb}{0.29,0.33,0.13}
\definecolor{um0color}{rgb}{1.1,0.01,0.24}
\definecolor{umcolor}{rgb}{1.0,0.01,0.24}
\definecolor{nuscolor}{rgb}{0.0,0.3,0.62}

\newcommand\blue[1]{{\color{blue} #1}}

\usepackage{amssymb,pifont}

\usepackage{enumitem}
\setlist[enumerate]{leftmargin=.5in}
\setlist[itemize]{leftmargin=.5in}

\newlist{todolist}{itemize}{2}
\setlist[todolist]{label=$\square$}

\usepackage{tikz}
\usetikzlibrary{patterns}
\usetikzlibrary{shapes,calc,shapes,arrows}

\newcommand{\dr}{\mathrm{d} \mathbf{r}}

\title{A novel splitting strategy to accelerate solving generalized eigenvalue problem from Kohn--Sham density functional theory\tnoteref{t1}}
\tnotetext[t1]{Submitted date: \today}
\author[1]{Yang Kuang}
\ead{ykuang@gdut.edu.cn}
\author[2,3]{Guanghui Hu \corref{cor1}}
\ead{garyhu@um.edu.mo}
\cortext[cor1]{Corresponding author. Email:{\tt garyhu@um.edu.mo}}

\address[1]{School of Mathematics and Statistics \&  Center for Mathematics and Interdisciplinary Science (CMIS), Guangdong University of Technology,  China}
\address[2]{State Key Laboratory of Internet of Things for Smart City and Department of Mathematics, University of Macau, Macao SAR, China.}
\address[3]{Zhuhai UM Science and Technology Research Institute, Zhuhai, China.}

\begin{document}

\begin{abstract}
In this paper, we propose a novel eigenpair-splitting method, inspired by the divide-and-conquer strategy, for solving the generalized eigenvalue problem arising from the Kohn-Sham equation. Unlike the commonly used domain decomposition approach in divide-and-conquer, which solves the problem on a series of subdomains, our eigenpair-splitting method focuses on solving a series of subequations defined on the entire domain. This method is realized through the integration of two key techniques: a multi-mesh technique for generating approximate spaces for the subequations, and a soft-locking technique that allows for the independent solution of eigenpairs. Numerical experiments show that the proposed eigenpair-splitting method can dramatically enhance simulation efficiency, and its potential towards practical applications is also demonstrated well through an example of the HOMO-LUMO gap calculation. Furthermore, the optimal strategy for grouping eigenpairs is discussed, and the possible improvements to the proposed method are also outlined.

\end{abstract}
\begin{keyword}
	Kohn--Sham equation; eigenpair-splitting; multi-mesh adaptive method; soft-locking strategy; 
\end{keyword}

\maketitle

\section{Introduction}
The Kohn--Sham density functional theory plays an important role within the realms of quantum physics, condensed matter, and computational chemistry \cite{kohn1965self}. It provides an effective approach to determine the ground state of a  quantum system containing $N_{nuc}$ nuclei and $N_{ele}$ electrons by solving the $N_{occ}$ lowest eigenpairs ${(\varepsilon_i,\psi_i)},{i=1,\dots,N_{occ}}$ from the Kohn--Sham (KS) equation, which is given by
\begin{equation}\label{eq:ks}
	\left(-\frac{1}{2}\nabla^2+V_\mathrm{ext}+V_\mathrm{Har}+V_\mathrm{xc} \right)\psi_i(\mathbf{r}) = \varepsilon_i\psi_i(\mathbf{r}).
\end{equation}
In this equation, $N_{occ}$ is defined as $N_{ele}/2$ for even $N_{ele}$ and $(N_{ele}+1)/2$ for odd, representing the number of occupation orbitals. The external potential  $V_\mathrm{ext}(\mathbf{r})$ is characterized by the sum of the Coulombic interactions between the electrons and the nuclei positioned at $\{\mathbf{R}_I\}_{I=1,\dots,N_{nuc}}$ with  respective charges $\{Z_I\}_{I=1,\dots,N_{nuc}}$, expressed as 
\begin{equation*}
	V_\mathrm{ext}(\mathbf{r})= -\sum_{I=1}^{N_{nuc}} \frac{Z_I}{|\mathbf{R}_I-\mathbf{r}|}.
\end{equation*}
The Hartree potential term  $V_\mathrm{Har}(\mathbf{r})$  encapsulates the classical electrostatic repulsion among the electrons. Lastly, $V_\mathrm{xc}(\mathbf{r})$ denotes the exchange-correlation potential, which accounts for the quantum mechanical effects of exchange and correlation that distinguish the Kohn-Sham non-interacting fictitious system from the actual many-electron system. Due to the lack of the analytical solution for the KS equation, it becomes necessary to employ numerical methods to obtain the solution of this equation. 

There exists several challenges in solving the KS equation. Firstly, the singularities in external potential term must be carefully addressed. The pseudopotential approaches \cite{pickett1989pseudopotential} provide a way to circumvent the singularities by dividing electrons into valence electrons and inner core electrons, thereby reducing the atom to an inoic core that interacts with the valence electrons. However, under extreme environments, pseudopotential approaches may lead to mispredictions \cite{xiao2010first, hamann2013optimized}.  Consequently, all-electron calculations that treat the external potential exactly are often necessary. To handle the singularities, the adaptive finite element (AFE) methods which enable different mesh sizes are favored \cite{ainsworth1997posteriori, verfurth2013posteriori}. Additionally, AFE methods have become increasingly attractive for electronic structure calculations in recent decades due to their ability to manage non-periodic boundary conditions and complex computational domains \cite{tsuchida1995electronic, lehtovaara2009all, bao2012h, fang2012kohn, motamarri2013higher,chen2014adaptive, maday2014h, davydov2016adaptive, motamarri2020dft}. Secondly, the KS equation is essentially an eigenvalue problem. As the size of the quantum system increases, the diagonalization and orthogonalization of the eigenvalue problem become computationally expensive \cite{lin2019numerical}. Specifically, the size of the discretized system and the number of required eigenpairs both increase, making the calculation more challenging.

Many efforts have been made to reduce the computational cost, such as orthogonalization-free methods \cite{dai2020gradient, gao2022orthogonalization} and the divide-and-conquer (DAC) method \cite{yang1991direct1, yang1991direct2, chen2016analysis}. In the latter approach, the global system is divided into several subsystems, and each subsystem is solved separately. This involves the decomposition of the computational domain. Inspired by this, we propose an alternative approach to divide the global system. Instead of decomposing the domain, we split the problem based on the eigenpairs. Specifically, the KS equation is decomposed into several subequations, each associated with a group of the eigenpairs. These subequations are solved independently to obtain the corresponding eigenpairs. We refer to this approach as the \emph{splitting method}.


The splitting of the eigenpairs based on the difference of the regularity among the corresponding wavefunctions, similar to the idea behind pseudopotentials. Specifically, the wavefunctions of inner core electrons vary more rapidly than those of valence electrons. Moreover, the decay of valence wavefunctions is slower than that of inner core wavefunctions. Therefore, when applying the adaptive finite element discretization to the KS equation, the mesh should effectively capture all the wavefunctions. Although the finite element space constructed on such a mesh is able to capture the variations for all wavefunctions, it may not be the most optimal space for any certain wavefunction.  Based on this observation, we aim to tailor the discretized space for each group of wavefunctions which share similar regularity, meaning the errors of these waveftunctions are similar. This idea can be implemented using the multi-mesh adaptive strategy \cite{li2005multi, kuang2024towards}.

Unlike the method that solves all equations on the same mesh (hereafter referred to as the \textit{single-mesh method}), the multi-mesh adaptive method (abbreviated as the \textit{multi-mesh method}) solves different variables with varying regularities on different meshes within the same computational domain. The mesh quality for a specific variable is improved through a mesh adaptation process similar to that used in the single-mesh method. In this study, the multi-mesh method is implemented based on the framework proposed in \cite{li2005multi}. Within this framework, mesh adaptation is achieved by locally refining or coarsening elements, a process known as $h$-adaptation. Specifically, \emph{a posteriori} error estimation is employed to guide the adaptations for all the meshes. The error indicator for each mesh is derived from the respective eigenpairs. It is noteworthy that the multi-mesh method has been used to handle KS wavefunctions and the Hartree potential \cite{kuang2024towards}, but this work differs in its approach to managing the KS wavefunctions.

Several numerical challenges remain in the splitting method for the KS equation, including the splitting strategy and the solution of eigenvalue problems on different spaces. While one could further divide the eigenpairs into more groups beyond just valence and inner core electrons, we found that further splitting can negatively impact the convergence of the self-consistent iteration and increase storage requirements. For solving the eigenvalue problems, we indeed need to solve two separate problems. In this work, we use the LOBPCG method \cite{knyazev2001toward} to solve the discretized eigenvalue problem. For the first problem, we use traditional methods to find the smallest eigenpairs, while for the second problem, we need to solve for the middle eigenvalues, which is more complex. We employ the soft-locking strategy \cite{knyazev2007block} to address this issue.

Additionally, the orthogonality of the wavefunctions must be carefully maintained. Since the wavefunctions belong to different finite spaces, orthogonality can be compromised, which often occurs in simulations. The loss of orthogonality primarily arises from the interpolation and projection processes. However, the resulting physical quantities, such as the total energy, are not significantly affected by this loss, as verified by several examples. To further ensure the validity of our results, we present a technique to guarantee orthogonality at a negligible cost.

In the next section, we introduce the KS model and the  finite element discretizations. Section 3 details the splitting method and the associated numerical challenges. Section 4 presents various numerical examples. Finally, we conclude the paper with a summary of our work.

\section{An $h$-adaptive finite element framework for Kohn--Sham equation}
\subsection{Kohn--Sham equation}
We rewrite the Kohn--Sham equation as the following eigenvalue problem
\begin{equation}\label{eq:KS}
	\left\{ 
	\begin{array}{lr}
		\hat{H}\psi_l(\mathbf{r}) = \varepsilon_l\psi_l(\mathbf{r}),& l=1,2,\dots , N_{occ},\\
		\displaystyle\int_{\mathbb{R}^3} \psi_l\psi_{l'} \,\mathrm{d}\mathbf{r} = \delta_{ll'},
		&l,l'=1,2,\dots ,N_{occ}, 
	\end{array} \right. 
\end{equation}
where the Hamiltonian contains the kinetic operator and potential terms
\begin{equation}\label{eq:Hamiltonian}
	\hat{H}([\rho];\mathbf{r}) = -\frac{1}{2}\nabla_{\mathbf{r}}^2 + V_{\mathrm{ext}}(\mathbf{r}) + V_{\mathrm{Har}}([\rho];\mathbf{r}) + V_{\mathrm{xc}}([\rho];\mathbf{r}).
\end{equation}
Here we introduce the notation $V([\rho];\mathbf{r})$ which implies that $V$ is a functional of the electron density $\rho(\mathbf{r}) = \sum_{l=1}^{N_{occ}}f_i|\psi_l|^2(\mathbf{r})$ with $f_i$ the occupation number for the $i$-th orbital.

The third term in $\hat{H}$  is the Hartree potential describing the Coulomb repulsion among the electrons
\begin{equation}\label{eq:har}
	V_{\mathrm{Har}}([\rho];\mathbf{r})=\int_{\mathbb{R}^{3}} \frac{\rho(\mathbf{r'})}
	{\lvert \mathbf{r}-\mathbf{r'} \lvert} \,\mathrm{d}\mathbf{r'}.
\end{equation}
The last term $V_{\mathrm{xc}}$ stands for the exchange-correlation potential, which is caused by the Pauli exclusion principle and other non-classical Coulomb interactions. Note that the analytical expression for the exchange-correlation term is unknown and therefore an approximation is needed. Specifically, the local density approximation (LDA) from the library Libxc \cite{marques2012libxc} with the slater exchange potential and the Vosko-Wilk-Nusair (VWN4) \cite{vosko1980accurate} is adopted in this work.  

Note that direct evaluation of the Hartree potential \eqref{eq:har} requires computational cost $\mathcal{O}(N^2)$ with $N$ being the number of grid points on the computational domain $\Omega$. For simplicity, we denote $\phi = V_{\mathrm{Har}}(\mathbf{r})$ hereafter. In this paper, the Hartree potential is obtained by solving the Poisson equation
\begin{equation}\label{eq:poisson}
	\left\{
	\begin{aligned}
		-\nabla^2\phi(\mathbf{r}) = 4\pi \rho(\mathbf{r}),&~ \text{in } \Omega,\\
		\phi(\mathbf{r}) = \phi_{\partial \Omega}(\mathbf{r}),&~ \text{on } \partial \Omega.
	\end{aligned}
	\right.
\end{equation}
The boundary value $\phi_{\partial \Omega}(\mathbf{r})$ is evaluated by the multipole expansion method. Specifically, the following approximation is used 
\begin{equation*}
	\begin{aligned}
		\left.
		\phi(\mathbf{r})\right|_{\mathbf{r} \in \partial \Omega} \approx  \frac{1}{\left|\mathbf{r}-\mathbf{r}^{\prime \prime}\right|} \int_{\Omega} \rho\left(\mathbf{r}^{\prime}\right) \,\mathrm{d} \mathbf{r}^{\prime}+\sum_{i=1,2,3} p_{i} \cdot \frac{r^{i}-r^{\prime \prime, i}}{\left|\mathbf{r}-\mathbf{r}^{\prime \prime}\right|^{3}} \\+\sum_{i, j=1,2,3} q_{i j} \cdot \frac{3\left(r^{i}-r^{\prime \prime, i}\right)\left(r^{j}-r^{\prime \prime, j}\right)-\delta_{i j}\left|\mathbf{r}-\mathbf{r}^{\prime \prime}\right|^{2}}{\left|\mathbf{r}-\mathbf{r}^{\prime \prime}\right|^{5}},
	\end{aligned}
\end{equation*}
where
\begin{equation*}
	\begin{aligned}		
	p_{i}&=\int_{\Omega} \rho\left(\mathbf{r}^{\prime}\right)\left(r^{\prime, i}-r^{\prime \prime, i}\right) \,\mathrm{d} \mathbf{r}^{\prime}, \\
	\quad q_{i j}&=\int_{\Omega} \frac{1}{2} \rho\left(\mathbf{r}^{\prime}\right)\left(r^{\prime, i}-r^{\prime \prime, i}\right)\left(r^{\prime, j}-r^{\prime \prime, j}\right) \,\mathrm{d} \mathbf{r}^{\prime}.
	\end{aligned}
\end{equation*}
In the above expressions, $\mathbf{r}^{\prime\prime}$ stands for an arbitrary point in $\Omega$. In the simulations, we choose it to be
\begin{equation*}
	\mathbf{r}^{\prime\prime} = \frac{\int \mathbf{r} \rho(\mathbf{r}) \,d\mathbf{r}}{\int  \rho(\mathbf{r}) \,d\mathbf{r}}.
\end{equation*}
As a result, $p_i=0, i = 1,2,3$, and the boundary condition can be simplified.

\subsection{Finite element discretization}
In practical simulations, a bounded polyhedral domain $\Omega \subset \mathbb{R}^3$ is served as the computational domain. Assume that the finite element space $V_h$ is constructed on $\Omega$ and the finite element basis of $V_h$ is denoted as $\{\varphi_1,\dots,\varphi_n\}$ with $n$ being the dimension of the space. Then the wavefunction $\psi_i$ can be approximated as $\psi_i^h$ on $V_h$ via
\begin{equation}\label{eq:x-fem}
	\psi_{i}^{h}(\mathbf{r})=\sum_{k=1}^{n} X_{k, i} \varphi_{k}, \quad X_{k, i}=\psi_{i}\left(\mathbf{r}_{k}\right), \quad X \in \mathbb{R}^{n \times p},
\end{equation}
where $\mathbf{r}_k$ denotes the node corresponding to the $i$-th basis function. As a result, to find the approximation of the wavefunctions $\{\psi_i\}$ on $V_h$ is to find $X$, i.e., the value of $\psi_i$ at each node $\mathbf{r}_k$.

In the finite element space $V_h$, the variation form for the KS equation \eqref{eq:KS} becomes: find $(\psi_i^h,\varepsilon_i^h)_{i=1,\dots,p} \in V_{h}\times \mathbb{R}$, such that
\begin{equation*}
	\int_{\Omega}\varphi\hat{H}\psi_i^h
	\,\mathrm{d} \mathbf{r} = \varepsilon_i^h\int_{\Omega} \varphi \psi_i^h \,\mathrm{d} \mathbf{r},\quad \forall \varphi \in V_h.
\end{equation*}
By letting $\varphi$ be the finite element basis function and inserting \eqref{eq:x-fem} to the above variational form, we have the following discrete eigenvalue problem
\begin{equation}\label{eq:evp-vh}
	A X = \varepsilon M Y.
\end{equation}
Here $A$ and $M$ are symmetric matrices with the entries
\begin{align}
	\label{eq:matA}
	A_{i,j} &=  \int_{\Omega} \frac{1}{2}\nabla \varphi_j\cdot \nabla \varphi_i + \big( V_{\mathrm{ext}}+\phi+V_{\mathrm{xc}} \big)\varphi_j\varphi_i \,\mathrm{d} \mathbf{r},\\
	\label{eq:matM}
	M_{i,j} &= \int_{\Omega} \varphi_j\varphi_i \,\mathrm{d} \mathbf{r}.
\end{align}

Similarly, we can obtain the linear system for the Poisson equation \eqref{eq:poisson} on $V_h$:
\begin{equation}\label{eq:ls-har}
	S \Phi = \mathbf{b},
\end{equation}
where $S$ is the stiff matrix with the entry
\begin{equation}\label{eq:matS}
	S_{i,j} =  \int_{\Omega} \frac{1}{2}\nabla \varphi_j\cdot \nabla \varphi_i \,\mathrm{d} \mathbf{r},
\end{equation}
and the right hand side $\mathbf{b}$ and the discretized $\phi^h$ 
\begin{equation*}
	\mathbf{b}_i =  \int_{\Omega} 4\pi \rho(\mathbf{r}) \varphi_i \,\mathrm{d} \mathbf{r},\quad
	\phi^h(\mathbf{r}) = \sum_{k=1}^n \Phi_k \varphi_k,\quad \Phi \in \mathbb{R}^n.
\end{equation*}

Owing to the singularity in the Hamiltonian \eqref{eq:Hamiltonian}, a uniform finite element discretization would lead to a large number of mesh grids to achieve chemical accuracy. Hence the adaptive mesh method is necessary to efficiently solve the Kohn--Sham equation, which will be discussed in the next section.

\subsection{Adaptive finite element method}
The adaptive mesh techniques offer enhanced numerical accuracy while requiring fewer mesh grids compared to uniform mesh methods. In this study, our primary focus lies on the $h$-adaptive methods, which involve local refinement and/or coarsening of mesh grids. An important aspect of $h$-adaptive methods is the determination of an error indicator. Generally speaking, error indicators identify regions within the domain that necessitate local refinement or coarsening, and they are typically derived from \emph{a posteriori} error estimations \cite{verfurth2013posteriori}. When there is only one orbital in the system, it is natural to generate the indicator based on information from that specific orbital. However, when there are multiple orbitals in the system, generating the indicator solely from an individual orbital is no longer advisable. This is because every orbital in the system is expected to be well-resolved using the mesh grids after mesh adaptation. To address this, we adopt the strategy proposed in \cite{bao2012h} for indicator generation. First, indicators are individually generated for each orbital using a specific method. Then, normalization is applied to each indicator. The final indicator is obtained by combining these normalized indicators.

Specifically, based on the \emph{a posteriori} error estimates \cite{verfurth2013posteriori},  the residual-based \emph{a posterirori} error indicator  for the KS equation \eqref{eq:KS} in the element $\mathcal{T}_K$ can be defined as
\begin{equation}\label{eq:indi-ks}
	\begin{aligned}
	\eta_{K,\mathrm{KS}} =\Big(&h_{K}^{2}\sum_{l=1}^p\big\|{\color{black}\mathcal{R}}_{K,\mathrm{KS}}(\psi_l)\big\|_{K}^{2}\\
	&+\sum_{e \subset \partial \mathcal{T}_{K}} \frac{1}{2} h_{e}\sum_{l=1}^p\big\|\mathcal{J}_{e}(\psi_l)\big\|_{e}^{2}\Big)^{\frac{1}{2}},	
	\end{aligned}
\end{equation}
where $h_K$ represents the largest length of the edges of the element $\mathcal{T}_K$, $h_{e}$ stands for the largest length of the common face $e$ of $\mathcal{T}_K$ and $\mathcal{T}_J$, ${\color{black}\mathcal{R}}_{K,\mathrm{KS}}(\psi)$ and $\mathcal{J}_{e}(\psi)$ are the residual and jump term on the element $\mathcal{T}_K$, whose formulations are written as
\begin{equation*}
	\left\{
	\begin{aligned}
		&{\color{black}\mathcal{R}}_{K,\mathrm{KS}}(\psi_l)=
		\hat{H} \psi_l- \varepsilon_lM\psi_l,\\
		&\mathcal{J}_{e}(\psi_l) = (\nabla
		\psi_l\bigm|_{\mathcal{T}_{K}} - \nabla\psi_l
		\bigm|_{\mathcal{T}_{J}}) \cdot \mathbf{n}_{e},
	\end{aligned}
	\right.
\end{equation*}
where $\mathbf{n}_{e}$ is the out normal vector on the face $e$ w.r.t the element $\mathcal{T}_K$. The definition of the norms are
\begin{align*}
	\left\| f(x) \right\|_K = \left(\int_K (f(x))^2dx\right)^{\frac{1}{2}},\\
	\left\| f(x) \right\|_{e} = \left(\int_{e} (f(x))^2dx\right)^{\frac{1}{2}}.
\end{align*}
The indicator \eqref{eq:indi-ks} involves the error arising from the Kohn--Sham equation \eqref{eq:KS} and is usually directly adopted to guide the mesh adaption. 

An aspect that the indicator \eqref{eq:indi-ks} possibly overlooks is the error associated with the Hartree potential. It is important to emphasize the accurate approximation of the Hartree potential, as it is a crucial component of the Hamiltonian \eqref{eq:Hamiltonian} and contributes to the Hartree potential energy. As the Hartree potential is obtained by solving the Poisson equation \eqref{eq:poisson}, a similar approach can be employed to generate an error indicator specifically for the Hartree potential, just as done for the Kohn--Sham equation. Specifically, the error indicator for the Hartree potential can be defined as
\begin{equation}\label{eq:indi-poisson}
	\eta_{K,\mathrm{Har}} =\left(h_{K}^{2}\big\|{\color{black}\mathcal{R}}_{K,\mathrm{Har}}(\phi)\big\|_{K}^{2}+\sum_{e \subset \partial \mathcal{T}_{K}} \frac{1}{2} h_{e}\big\|\mathcal{J}_{e}(\phi)\big\|_{e}^{2}\right)^{\frac{1}{2}},
\end{equation}
where the residual part is defined as ${\color{black}\mathcal{R}}_{K,\mathrm{Har}}(\phi)= \nabla^2\phi + 4\pi\rho$. Based on the  analysis above,  the second indicator for solving the Kohn--Sham equation can be designed as 
\begin{equation} \label{eq:indi-kshar}
	\eta_{K,\mathrm{KS+Har}} = \sqrt{\eta_{K,\mathrm{Har}}^2 + \eta_{K,\mathrm{KS}}^2}.
\end{equation}
The aim of using this error indicator is to generate a mesh on which both the wavefunctions and the Hartree potential are approximated well. 

The mesh adaptation process, along with the solution of the Kohn--Sham equation, can be carried out using either the first \eqref{eq:indi-ks} or the second error indicator \eqref{eq:indi-kshar}. It is important to note that the adaptive algorithms utilizing these two different indicators are essentially identical, except for the choice of the indicator. Moreover, all simulations are performed on a single mesh simultaneously, i.e., the Poisson equation \eqref{eq:poisson} for the Hartree potential is discretized and solved on the same finite element space for the Kohn--Sham equation.

\section{A novel splitting strategy for the calculation of Kohn-Sham orbitals}
\subsection{The motivation and idea for the orbital splitting}\label{subsec:split}
\begin{figure*}[!htp]
	\centering
	\includegraphics[width=0.32\linewidth]{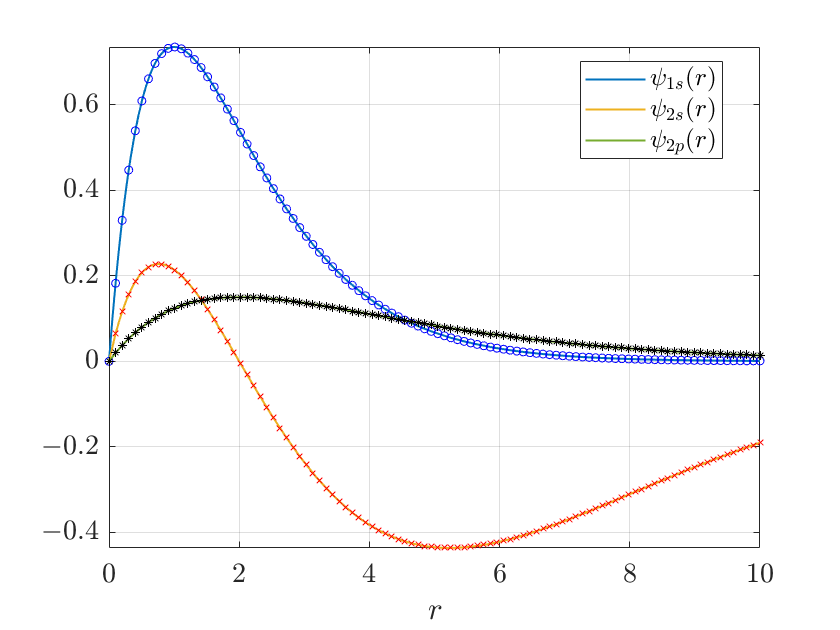}
	\includegraphics[width=0.32\linewidth]{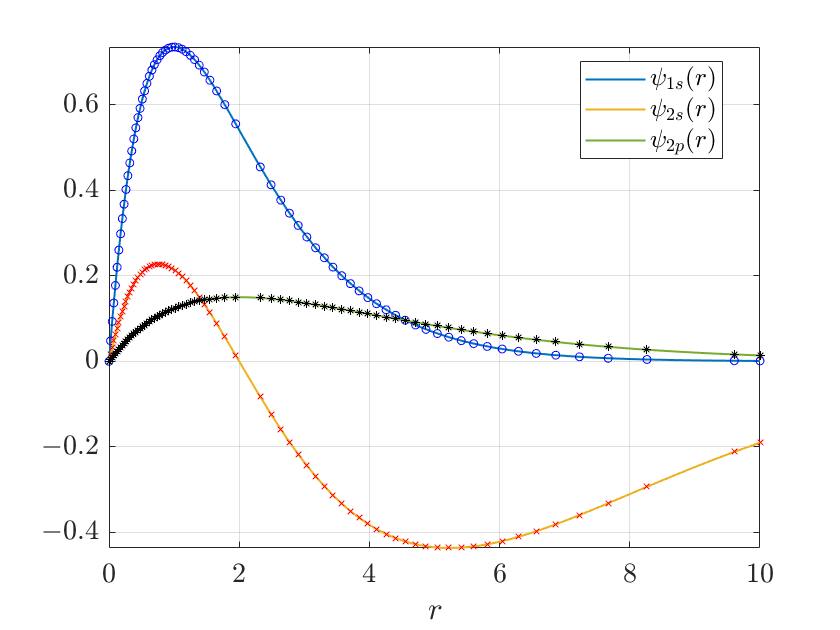}
	\includegraphics[width=0.32\linewidth]{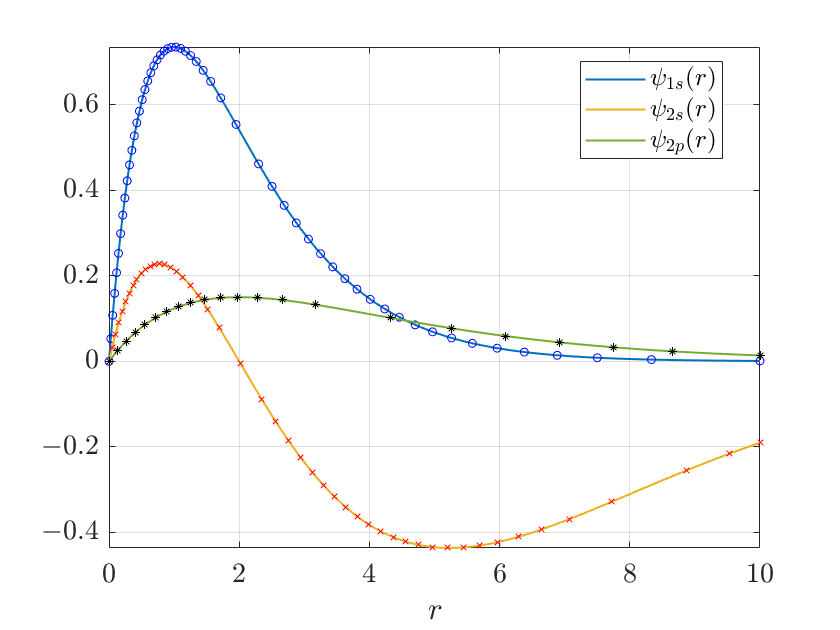}
	\caption{Three states for the hydrogen atom. \label{fig:hydrogen}}
\end{figure*}

For a clear illustration of our idea, we start with considering the linear finite element approach to hydrogen atom. Its governing equation is 
\begin{equation}
	\left(-\frac{1}{2}\nabla^2 - \frac{1}{r}\right)\psi_i  = \varepsilon_i \psi_{i}
\end{equation}
with first three normalized radial solutions \cite{froese1997computational} as
\begin{align*}
	&\psi_{1s} = 2re^{-r},
	~\psi_{2s} = \frac{r}{2}e^{-r/2}(1-\frac{1}{2}r),~\psi_{2p} = \frac{r^2}{2\sqrt{6}}e^{-r/2}.
\end{align*}
Now we ignore the solution part, and focus on the finite element approximation to these solutions. We test three different mesh strategies, i.e., the uniform mesh, the adaptive mesh which is based on equally distributed L2 error, and the multi-mesh adaptive method, which utilize a specific mesh to each orbital. The L2 error for the linear interpolation of the radial orbitals is presented in \Cref{tab:1dhydrogen}. The results demonstrate that to achieve an error within 0.003, the multi-mesh method outperforms the other two methods, which use fewer mesh grids. The corresponding wavefunctions and mesh grids are illustrated in \Cref{fig:hydrogen}. From this figure we find that the regularity of these three states are different. Specifically, the 1s state varies most rapid at the origin, while the other states vary slightly. Hence, to capture the variations for these three states, the mesh should be dense enough around the origin, as indicated in the middle figure of \Cref{fig:hydrogen}. Whereas, this mesh may not be the most suitable mesh for the other two states. This motivates us to use different meshes for the different states. As indicated in the right figure of \Cref{fig:hydrogen}, the multi-mesh method can achieve the same accuracy with less number of mesh grids.

\begin{table}
	\caption{L2 error for linear interpolation of the radial orbitals. $N_{1s},N_{2s},N_{2p}$ indicate the number of mesh grids for the respective orbitals. $\Delta$ represents the L2 error for each corresponding orbital. `uni-mesh' refers to the method using uniform mesh, `ada-mesh' denotes the method with an adaptively changed mesh, and `multi-mesh' means the method using a distinct mesh for each orbital. \label{tab:1dhydrogen}}
	\centering
	{\small
	\begin{tabular}{c|ll|ll|ll}
		\toprule
		&$N_{1s}$ &$\Delta_{1s}$ &$N_{2s}$&$\Delta_{2s}$&$N_{2p}$&$\Delta_{2p}$ \\ \midrule
		uni-mesh & 100 & 0.00208 & 100 & 0.00084 & 100 & 0.00015 \\
		ada-mesh &74 & 0.00059 & 74 & 0.00130 & 74 & 0.00063 \\
		mul-mesh & 54 & 0.00098 & 48 & 0.00091 & 21 & 0.00098 \\ 
		\bottomrule
	\end{tabular}
	}
\end{table}

The above atomic example demonstrates the idea of using the multi-mesh method. This idea can also be extended to the molecule cases. Before introducing this method, we first discuss the splitting of the orbitals, i.e., separate the orbitals into several groups. As is well known, to obtain the ground state from the Kohn--Sham equation, we have to calculate $p$ eigenpairs from the eigenvalue problem. Based on the observation of the distinct regularities of the wavefunctions, it is possible to separate the orbitals, similar to the approach used in the pseudopotential method.

In an atom, there are two types of electron orbitals: core electrons and valence electrons. The wavefunctions of core electrons exhibit strong variations near the nucleus and decay rapidly to zero far from the nucleus. In contrast, valence electron wavefunctions are relatively smooth near the nucleus and decay more slowly at greater distances. Therefore, it is straightforward to split the eigenpairs into these two groups for an atom.

This strategy can also be applied to molecules. We first collect the core electrons and then distribute the orbitals according to the principle that each orbital is occupied by two electrons. For example, in the water molecule (\ch{H2O}), there are 2 core electrons (from the oxygen atom) and 8 valence electrons, leading to a splitting of the eigenpairs into 1 core eigenpair and 4 valence eigenpairs.

An alternative way which only based on the numerical results is also available. Specifically, we solve the KS equation on an initial quality mesh for all the eigenpairs, and then we spit these eigenpairs based on the distribution of the eigenvalues.

We list several splitting examples here in \Cref{tab:split-mole}.
\begin{table}[!htp]
	\centering
	\caption{Splitting of the molecules \label{tab:split-mole}}
	\begin{tabular}{c|c|l|l|l}
		\toprule
		molecule &$N_\mathrm{ele}$&$p$& $p_1$ & $p_2$ \\ \midrule
		LiH & 4 & 2 & 1 &1 \\
		\ch{BeH2} & 6 & 3 & 1 & 2 \\
		\ch{H2O} & 10 & 5 & 1 & 4 \\
		\ch{C6H6} &42 &21 & 6 &15 \\
 		\bottomrule		
	\end{tabular}
\end{table}

\subsection{A multi-mesh technique for the implementation}
To solve the eigenstates separately, we introduce an approach known as the multi-mesh adaptive method \cite{li2005multi, cai2024AFEPack}. The primary objective of this method is to strike a balance between achieving chemical accuracy and managing computational costs effectively. By utilizing multiple meshes instead of a single mesh, the multi-mesh adaptive method offers greater flexibility in adapting the mesh resolution to capture variations for different quantities of interest. 
	
Specifically, in this method, we will utilize two distinct meshes for the eigenvalue problem. Both two meshes will be adapted during the simulation. Each mesh is specifically designed for solving the associated group eigenpairs of the Kohn--Sham equation. The mesh is tailored to ensure high resolution in regions where the respective group wavefunctions exhibit significant changes. By employing two separate meshes with specific focus areas, we can ensure that each group of eigenpairs is solved with the appropriate level of accuracy and capture the variations unique to the associated wavefunctions.
	
To implement the multi-mesh method, careful handling of two components is crucial. The first component involves effectively managing the mesh grids, allowing for flexible local refinement or coarsening as needed. This management mechanism should also facilitate efficient solution updates from the old mesh to the new mesh. The second component focuses on ensuring efficient and accurate communication between different meshes, particularly in the calculation of integrals. Efficient communication protocols play a vital role in maintaining consistency and accuracy across the various meshes. These requirements can be fulfilled by the hierarchical geometry tree data structure, which will be introduced in detail in the following.
	
\paragraph{The hierachical geometry}	A well-designed data structure for the mesh grids is needed for an effective management mechanism. In the presented algorithm, the hierarchical geometry tree (HGT) \cite{li2005multi, bao2012h} is utilized. 
	
Firstly, the mesh structure is described hierarchically, which means that the mesh information is given from the lowest dimension (0-D, the points) to the highest dimension (3-D, the tetrahedron) hierarchically. An element such as a point for 0-D, an edge for 1-D, a triangle for 2-D, a tetrahedron for 3-D is called a geometry. In the hierarchical description of a tetrahedron, all geometries have a belonging-to relationship. For example, if an edge is one of the edges of a triangle, this edge belongs to this triangle. With this hierarchical structure, the geometry information of a tetrahedron can be referred to flexibly, and the refinement and coarsening of a mesh can also be implemented efficiently. 

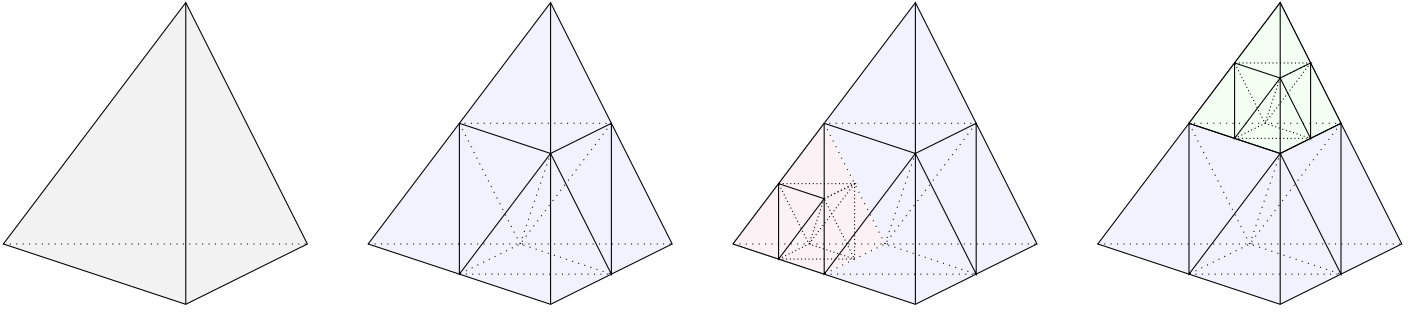
\begin{figure*}[h]
	\centering
	\begin{tikzpicture}
		\begin{scope}[scale=0.8]
			\draw[shift={(-6,0)},fill=gray!10] (0,0)--(3,-1)--(5,0)--(3,4)--cycle; 
			\draw[shift={(-6,0)}] (3,-1)--(3,4);
			\draw[shift={(-6,0)}, dotted] (0,0)--(5,0); 
			
			\draw[fill=blue!5] (0,0)--(3,-1)--(5,0)--(3,4)--cycle; 
			\draw[] (3,-1)--(3,4);
			\draw[dotted] (0,0)--(5,0); 
			
			\draw[] (1.5,-0.5)--(1.5,2)--(3,1.5)--cycle;
			\draw[] (3,1.5)--(4,2)--(4,-0.5)--cycle;
			\draw[dotted] (1.5,-0.5)--(4,-0.5)--(2.5,0)--cycle;
			\draw[dotted] (2.5,0)--(1.5,2)--(4,2)--cycle;
			\draw[dotted] (2.5,0)--(3,1.5);

			\draw[shift={(6,0)}, fill=blue!5] (0,0)--(3,-1)--(5,0)--(3,4)--cycle; 
			\draw[white,shift={(6,0)}, fill=purple!5] (0,0)--(1.5,-0.5)--(2.5,0)--(1.5,2)--cycle;
			\draw[shift={(6,0)}] (0,0)--(1.5,-0.5);
			\draw[shift={(6,0)}] (0,0)--(1.5,2);
			\draw[shift={(6,0)}] (3,-1)--(3,4);
			\draw[shift={(6,0)}, dotted] (0,0)--(5,0); 
			\draw[shift={(6,0)}] (1.5,-0.5)--(1.5,2)--(3,1.5)--cycle;
			\draw[shift={(6,0)}] (3,1.5)--(4,2)--(4,-0.5)--cycle;
			\draw[shift={(6,0)}, dotted] (1.5,-0.5)--(4,-0.5)--(2.5,0)--cycle;
			\draw[shift={(6,0)}, dotted] (2.5,0)--(1.5,2)--(4,2)--cycle;
			\draw[shift={(6,0)}, dotted] (2.5,0)--(3,1.5);
			\draw[shift={(6,0)}] (0.75,-0.25)--(1.5,0.75)--(0.75,1)--cycle;
			\draw[shift={(6,0)}, densely dotted] (1.5,0.75)--(2,1)--(2,-0.25)--cycle;
			\draw[shift={(6,0)}, densely dotted] (1.5,0.75)--(2,1)--(0.75,1);
			\draw[shift={(6,0)}, densely dotted] (2,1)--(1.25,0)--(0.75,1);
			\draw[shift={(6,0)}, densely dotted] (1.5,0.75)--(1.25,0);
			\draw[shift={(6,0)}, densely dotted] (0.75,-0.25)--(2,-0.25)--(1.25,0)--cycle;

			\draw[shift={(12,0)},fill=blue!5] (0,0)--(3,-1)--(5,0)--(3,4)--cycle; 
			\draw[shift={(13.5,2)}, fill=green!5] (0,0)--(1.5,-0.5)--(2.5,0)--(1.5,2)--cycle;

			\draw[shift={(12,0)}] (3,-1)--(3,4);
			\draw[shift={(12,0)}, dotted] (0,0)--(5,0); 
			\draw[shift={(12,0)}] (1.5,-0.5)--(1.5,2)--(3,1.5)--cycle;
			\draw[shift={(12,0)}] (3,1.5)--(4,2)--(4,-0.5)--cycle;
			\draw[shift={(12,0)}, dotted] (1.5,-0.5)--(4,-0.5)--(2.5,0)--cycle;
			\draw[shift={(12,0)}, dotted] (2.5,0)--(1.5,2)--(4,2)--cycle;
			\draw[shift={(12,0)}, dotted] (2.5,0)--(3,1.5);
			\draw[shift={(13.5,2)}] (0.75,-0.25)--(1.5,0.75)--(0.75,1)--cycle;
			\draw[shift={(13.5,2)}] (1.5,0.75)--(2,1)--(2,-0.25)--cycle;
			\draw[shift={(13.5,2)}, densely dotted] (1.5,0.75)--(2,1)--(0.75,1);
			\draw[shift={(13.5,2)}, densely dotted] (2,1)--(1.25,0)--(0.75,1);
			\draw[shift={(13.5,2)}, densely dotted] (1.5,0.75)--(1.25,0);
			\draw[shift={(13.5,2)}, densely dotted] (0.75,-0.25)--(2,-0.25)--(1.25,0)--cycle;
		\end{scope}
	\end{tikzpicture}
	
	\caption{First: Mesh on the root tetrahedron. Seond: Mesh on the global refinement of the tetrahedron. Third and forth: the local refinements of second mesh. \label{fig:tet}}
\end{figure*}

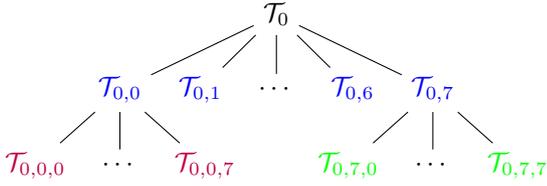
\begin{figure}[h]
	\centering
	\begin{forest}
		[$\mathcal{T}_{0}$
		[\blue{$\mathcal{T}_{0,0}$}
		[\color{purple}$\mathcal{T}_{0,0,0}$]
		[$\cdots$]
		[\color{purple}$\mathcal{T}_{0,0,7}$]
		]
		[\blue{$\mathcal{T}_{0,1}$}]
		[$\cdots$]
		[\blue{$\mathcal{T}_{0,6}$}]
		[\blue{$\mathcal{T}_{0,7}$}
		[\color{green}$\mathcal{T}_{0,7,0}$]
		[$\cdots$]
		[\color{green}$\mathcal{T}_{0,7,7}$]
		]
		]
	\end{forest}
	\caption{The octree data structure. \label{fig:octree}}
\end{figure}

Secondly, the mesh is stored and managed by a tree data structure. The validity of using the tree data structure is due to the strategy of element refinement and coarsening. Specifically, for a tetrahedron element (the left of \Cref{fig:tet}), the refinement of this tetrahedron is dividing it into eight equally small tetrahedrons via connecting the midpoints on each edge (the second column of \Cref{fig:tet}). As a result, a belonging-to relationship can be established, for example, any small tetrahedron that is called ``child" belongs to the large ``parent" tetrahedron. Meanwhile, the coarsening of any child tetrahedron in the second column of \Cref{fig:tet} is releasing all the children to obtain the parent tetrahedron. By this refinement and coarsening strategy, the tree data structure is able to be established, as indicated in \Cref{fig:octree}. 

The local refinements of a mesh can also be managed by the tree data structure. Assume the mesh is built on a root tetrahedron $\mathcal{T}_0$ (left of \Cref{fig:tet}), and the local refinement is implemented on the tetrahedron $\mathcal{T}_{0,0}$ (third column of \Cref{fig:tet}), then the tree data structure is shown in \Cref{fig:octree}. The root tetrahedron $\mathcal{T}_0$ has eight children, and the tetrahedron $\mathcal{T}_{0,0}$ has eight children $\{\mathcal{T}_{0,0,0},\cdots,\mathcal{T}_{0,0,7}\}$ as well. Remarkably, the mesh can be established from the octree data structure.  For example, the set $\{\mathcal{T}_0\}$ forms the mesh in the left of \Cref{fig:tet}, the set  $\{\mathcal{T}_{0,0},\mathcal{T}_{0,1},\mathcal{T}_{0,2},\mathcal{T}_{0,3}\}$ forms the second mesh in \Cref{fig:tet}, and the set $\{\mathcal{T}_{0,0,0},\mathcal{T}_{0,0,1},\cdots,\mathcal{T}_{0,0,7},\mathcal{T}_{0,1},\mathcal{T}_{0,2},\mathcal{T}_{0,3}\}$ forms the third mesh in \Cref{fig:tet}.

Such a data structure makes the management of multiple meshes available. Assume now for another variable, the local refinement of $\mathcal{T}_{0,0}$ is not necessary, on the contrary, the refinement on $\mathcal{T}_{0,7}$ is needed, then the mesh is shown in the right of \Cref{fig:tet}. In this case, the node of the octree $\mathcal{T}_{0,7}$ generates eight new children $\{\mathcal{T}_{0,7,0},\cdots,\mathcal{T}_{0,7,7}\}$, and this mesh is formed by the set $\{\mathcal{T}_{0,0},\mathcal{T}_{0,1},\cdots,\mathcal{T}_{0,6},\mathcal{T}_{0,7,0},\cdots,\mathcal{T}_{0,7,7}\}$, which can also be read from the octree \Cref{fig:octree}.

\begin{figure}[!h]
	\centering
	\begin{tikzpicture}
		\begin{scope}[scale=0.6]		
			\path[shift={(-7,0)}]  (3,4) coordinate (A) node [above] {$A$};
			\path[shift={(-7,0)}]  (0,0) coordinate (B) node [left] {$B$};
			\path[shift={(-7,0)}]  (5,0) coordinate (C) node [right] {$C$};
			\path[shift={(-7,0)}]  (3,-1) coordinate (D) node [below] {$D$};
			\path[shift={(-7,0)}]  (2.5,0) coordinate (E) node [above left] {$E$};
			\draw[shift={(-7,0)}] (0,0)--(3,-1)--(5,0)--(3,4)--cycle; 
			\draw[shift={(-7,0)}] (3,-1)--(3,4); 
			\draw[shift={(-7,0)}, densely dotted] (0,0)--(5,0); 
			\draw[shift={(-7,0)}, densely dotted] (3,4)--(2.5,0)--(3,-1);
			
			\path (3,4) coordinate (A) node [above] {$A$};
			\path (0,0) coordinate (B) node [left] {$B$};
			\path (5,0) coordinate (C) node [right] {$C$};
			\path (3,-1) coordinate (D) node [below] {$D$};
			\path (2.5,0) coordinate (E) node [above left] {$E$};
			\path (1.5,-0.5) coordinate (G) node [below left] {$G$};
			\path (4,-0.5) coordinate (F) node [below right] {$F$};
			\draw[] (0,0)--(3,-1)--(5,0)--(3,4)--cycle; 
			\draw[] (3,-1)--(3,4);
			\draw[densely dotted] (0,0)--(5,0); 
			\draw[] (3,4)--(1.5,-0.5);
			\draw[] (3,4)--(4,-0.5);
			\draw[densely dotted] (3,4)--(2.5,0); 
			\draw[densely dotted] (1.5,-0.5)--(4,-0.5)--(2.5,0)--cycle;
		\end{scope}
	\end{tikzpicture}
	\caption{Twin-tetrahedron and four-tetrahedron. \label{fig:twin-tet}}
\end{figure}
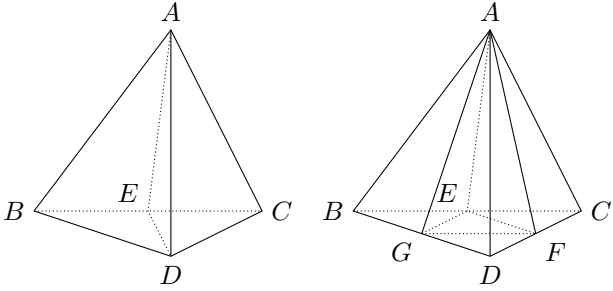

With the hierarchical description of the geometry and the tree data structure, meshes can be effectively managed by the HGT. However, building the finite element space directly on a mesh leads to non-conforming finite elements due to the hanging points in the direct neighbors of those refined tetrahedrons. To address these hanging points, two types of geometries are introduced: twin-tetrahedron and four-tetrahedron, as shown in \Cref{fig:twin-tet}. In the twin-tetrahedron geometry (\Cref{fig:twin-tet}, left), there are two tetrahedrons, ABED and AECD, with five degrees of freedom (DOF) at interpolation points A, B, C, D, and E. To ensure the finite element space is conforming, the basis functions are constructed as follows: each basis function has a value of 1 at its corresponding interpolation point and 0 at others. The support of the basis function for interpolation points common to both sub-tetrahedrons, such as A, D, and E, covers the entire twin-tetrahedron. For non-common points like B, the support is limited to tetrahedron ABED, while for point C, it is tetrahedron AECD. A similar strategy is applied to construct basis functions in the four-tetrahedron geometry (\Cref{fig:twin-tet}, right). By using these twin-tetrahedron and four-tetrahedron geometries, a conforming finite element space can be smoothly built in a mesh with local refinement.

In virtue of the tree data structure, multiple meshes are allowed to be described by the same HGT. Naturally, a problem arises: is it possible to make the information communicates among these meshes without any loss? The answer is affirmative, and the reason relies on the belonging-to relationship between any two nodes in the HGT. Briefly speaking, there exist only three kinds of relationship between any two nodes: equal, belonging-to, and no-overlap. The communication among equal or no-overlap elements is trivial. Consequently, we only need to take care of the second kind of relationship. We introduce the implementation details in the following.

\paragraph{The communications among different meshes}
In the presented method, the Hamiltonian matrices and energy require information from all the KS meshes and the Hartree mesh, which mainly involve numerical integration. We take the evaluation of the Hartree potential energy for example, which could be written as 
\begin{equation}\label{eq:har-cal}
	E_{\mathrm{Har}} = \frac{1}{2}\int_\Omega \phi(\mathbf{r}) \left(\sum_{l=1}^{p_1}\psi_l(\mathbf{r})^2 + \sum_{l=p_1+1}^p\psi_l(\mathbf{r})^2 \right)\dr. 
\end{equation}
The integration should be carefully treated since the Hartree potential  $\phi(\mathbf{r})$ and the wavefunction $\psi_l(\mathbf{r})$ belong to different finite element spaces built on different meshes. We first calculate the Hartree potential energy contributed from the first group. An intuitive illustration is given to show how the presented method calculates the first part of the integral \eqref{eq:har-cal}. For demonstration, we reuse the last two meshes in \Cref{fig:tet}. Assume the finite element space $V_{\mathcal{T}_{\mathrm{KS},1}}$ for the orbitals of the first group  is built on the mesh $\mathcal{T}_{\mathrm{KS},1}$ (the third column of \Cref{fig:tet}), and the space $V_{\mathcal{T}_{\mathrm{Har}}}$ is built on the mesh $V_{\mathcal{T}_{\mathrm{Har}}}$ (the last column of \Cref{fig:tet}). For the common elements such as $\mathcal{T}_{0,1}$ which belongs to both $\mathcal{T}_{\mathrm{KS},1}$ and $\mathcal{T}_{\mathrm{Har}}$,  the numerical integration can be calculated directly.  While for the remained elements, for example, in $\mathcal{T}_{\mathrm{KS},1}$, the tetrahedron $\mathcal{T}_{0,0}$ is refined, while in $\mathcal{T}_{\mathrm{Har}}$ this tetrahedron is kept, special treatment is needed to avoid the loss of accuracy. 

In order to prevent the loss of accuracy, we employ a strategy that maximizes the utilization of quadrature points in numerical integrals. Specifically, the numerical integration on the element $\mathcal{T}_{0,0}$ is divided into the integrations on its eight refined sub-tetrahedrons, 
\begin{align*}
&\frac{1}{2}\int_{\mathcal{T}_{0,0}} \phi(\mathbf{r}) \left(\sum_{l=1}^{p_1} \psi_l(\mathbf{r}) ^2 \right) \dr \\
&\approx \sum_{K=1}^8\sum_{j=1}^q \mbox{area}(K) J_j^K w_j^K \phi(\mathbf{r}_j^K)\left(\sum_{l=1}^{p_1} \psi_l(\mathbf{r}_j^K)^2 \right),
\end{align*}
where the element $K$ represent the eight sub-tetrahedrons of $\mathcal{T}_{0,0}$, $\mathbf{r}_j^K$ is the $j$-th quadrature point of $K$, $J_j^K$ is the jacobian at $\mathbf{r}_j^K$, and $w_l^K$ is the associated weight in the numerical quadrature. The values of the Hartree potential on these quadrature points can be obtained by numerical interpolation. Similarly, the integral on the tetrahedron $\mathcal{T}_{0,7}$ is evaluated by doing the numerical integral on its sub-triangles. In this way, the accuracy of the integral will not be affected in the communication among $\mathcal{T}_{\mathrm{KS},1}$ and $\mathcal{T}_{\mathrm{Har}}$. The second part of Hartree potential energy \Cref{eq:har-cal} can be evaluated similarly. 

With the HGT, the solution update from the old mesh to the new mesh after the local refinement and coarsening can also be implemented efficiently. As we mentioned before, a set of all leaf nodes forms a mesh. Although different meshes correspond to different sets of leaf nodes, all meshes are from the same set of root nodes. Then the belonging-to relationship of the geometries between two meshes can be analyzed easily. As a result, the solution update can be implemented efficiently according to the relationship.

\paragraph{Residual \emph{a posteriori} error estimate in multi-mesh method}
Assume the orbtials are splitted into two groups as described in \Cref{subsec:split}, with the corresponding wavefunctions denoted as $\{\psi_1,\cdots,\psi_{p_1}\}$ and $\{\psi_{p_1+1},\cdots,\psi_{p}\}$, respectively. The first group of wavefunctions is solved on the mesh $\mathcal{T}_1$, while the second group is solved on the mesh $\mathcal{T}_2$. Based on the \emph{a posteriori} error estimation \Cref{eq:indi-ks} for mesh adaption,  the error indicator for the first mesh can be designed as 
\begin{equation}\label{eq:indi-ks-1}
	\begin{aligned}
	\eta_{K,\mathrm{KS},1} =\Big(&h_{K}^{2}\sum_{l=1}^{p_1}\big\|{\color{black}\mathcal{R}}_{K,\mathrm{KS}}(\psi_l)\big\|_{K}^{2}\\
	&+\sum_{e \subset \partial \mathcal{T}_{K}} \frac{1}{2} h_{e}\sum_{l=1}^{p_1}\big\|\mathcal{J}_{e}(\psi_l)\big\|_{e}^{2}\Big)^{\frac{1}{2}},	
	\end{aligned}
\end{equation}
and the error indicator for the second mesh can be designed as
\begin{equation}\label{eq:indi-ks-2}
	\begin{aligned}
	\eta_{K,\mathrm{KS},2} =\Big(&h_{K}^{2}\sum_{l=p_1+1}^{p}\big\|{\color{black}\mathcal{R}}_{K,\mathrm{KS}}(\psi_l)\big\|_{K}^{2}\\
	&+\sum_{e \subset \partial \mathcal{T}_{K}} \frac{1}{2} h_{e}\sum_{l=p_1+1}^{p}\big\|\mathcal{J}_{e}(\psi_l)\big\|_{e}^{2}\Big)^{\frac{1}{2}}.
	\end{aligned}
\end{equation}
The only difference between these two error indicators is the wavefunctions they consider, as each indicator only accounts for the wavefunctions on its mesh.

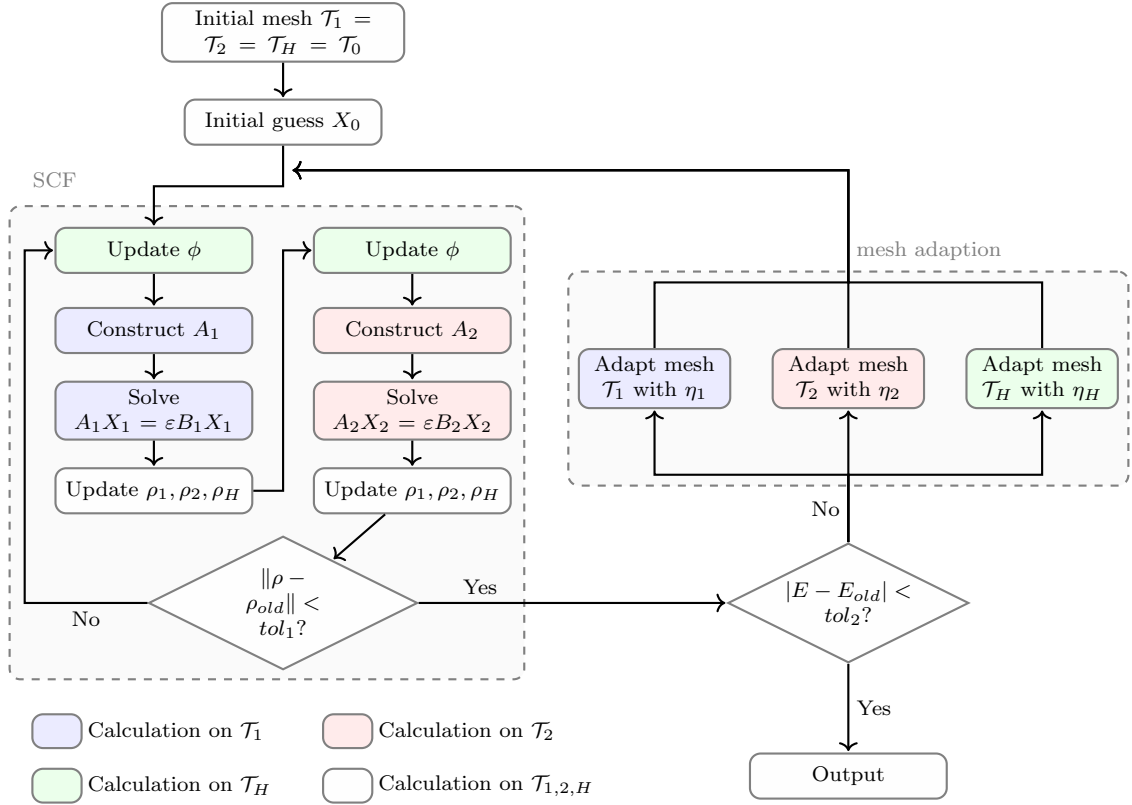
\begin{figure*}[htbp]
	\centering \footnotesize
	\begin{tikzpicture}[
		scale=.85,
		auto,
		decision/.style = { diamond, aspect=2, draw=gray,
			thick, fill=gray!1, text width=4em, text badly centered,
			inner sep=1pt},
		block/.style = { rectangle, draw=gray, thick, fill=gray!1,
			text width=8em, text centered, rounded corners,
			minimum height=2em },
		blockhar/.style = { rectangle, draw=gray, thick, fill=green!8,
			text width=8em, text centered, rounded corners,
			minimum height=2em },
		blockks1/.style = { rectangle, draw=gray, thick, fill=blue!8,
			text width=8em, text centered, rounded corners,
			minimum height=2em },
		blockks2/.style = { rectangle, draw=gray, thick, fill=red!8,
			text width=8em, text centered, rounded corners,
			minimum height=2em },
		line/.style = { draw, thick, ->, shorten >= 0.5pt},
		]
		
		\node [block, text width=10em] at (0,1.9) (initmesh) {Initial mesh $\mathcal{T}_1 = \mathcal{T}_2={\mathcal{T}}_H = \mathcal{T}_0 $};
		\node [block, text width=8em] at (0,0.5) (prob) {Initial guess $X_0$};

		\node at (0,-0.25) (null1) {};
		
		\node [blockks1, text width=1.5em,minimum height=1.5em] at (-3.5,-9) (ksmesh1) {};
		\node [right] at (-3.15,-9) {Calculation on $\mathcal{T}_1$};
		\node [blockks2, text width=1.5em,minimum height=1.5em] at (1.,-9) (ksmesh2) {};
		\node [right] at (1.35,-9) {Calculation on $\mathcal{T}_2$};
		\node [blockhar, text width=1.5em,minimum height=1.5em] at (-3.5,-9.85) (harmesh) {};
		\node [right] at (-3.15,-9.85) {Calculation on $\mathcal{T}_H$ };
		\node [block, text width=1.5em,minimum height=1.5em] at (1.,-9.85) (allmesh) {};
		\node [right] at (1.35,-9.85) {Calculation on $\mathcal{T}_{1,2,H}$};
		
		\node [block,rectangle, draw=gray, thick,dashed, fill=gray!4,
		text width=22em, text centered, rounded corners, minimum
		height=21em] at (-0.25,-4.5) (scf) {};
		\node [right, color=gray] at (-4,-0.4) {SCF};
		
		\node [block,rectangle, draw=gray, thick,dashed, fill=gray!4,
		text width=24em, text centered, rounded corners, minimum
		height=9.5em] at (8.75,-3.5) (ada) {};
		\node [right,color=gray] at (8.75,-1.5) {mesh adaption};

		\node [blockhar] at (-2,-1.5) (genhar) {Update $\phi$};
		\node [blockhar] at (2,-1.5) (genhar2) {Update $\phi$};
		\node [blockks1] at (-2,-2.75) (genh) {Construct ${A}_1$};
		\node [blockks1, text width = 8em] at (-2,-4)  (eigensolver) {Solve ${A}_1{X}_1 = {\varepsilon} {B}_1{X}_1$};
		\node [blockks2] at (2,-2.75) (genh2) {Construct ${A}_2$};
		\node [blockks2, text width = 8em] at (2,-4)  (eigensolver2) {Solve ${A}_2{X}_2 = {\varepsilon}{B}_2{X}_2$};
		\node [block] at (-2, -5.25) (updaterho) {Update $\rho_1,\rho_2,\rho_H$};
		\node [block] at (2, -5.25) (updaterho2) {Update $\rho_1,\rho_2,\rho_H$};
		\node [decision, text width=5em] at (0,-7) (scfdec) {$\|\rho-\rho_{old}\| < tol_1 ?$};
		
		\node [decision, text width=6em] at (8.75,-7) (adadec) {$|E-E_{old}| < tol_2 ?$};
		
		\node [blockks1,text width=6em] at (5.75,-3.5) (adaks1) {Adapt mesh $\mathcal{T}_1$ with $\eta_1$};
		
		\node [blockks2,text width=6em] at (8.75,-3.5) (adaks2) {Adapt mesh $\mathcal{T}_2$ with $\eta_2$};
		
		\node [blockhar,text width=6em] at (11.75,-3.5) (adahar) {Adapt mesh $\mathcal{T}_H$ with $\eta_H$};

		\node [block] at (8.75, -9.7) (end) {Output};


		\begin{scope} [every path/.style=line,thick,shorten >= 0.5pt]
			\path (initmesh) -- (prob);
			\path (prob)    --(0,-0.5)--(-2,-0.5)--   (genhar);   
			\path (genhar)    --  (genh);
			\path (genh) -- (eigensolver);
			\path (eigensolver) -- (updaterho);
			\path (updaterho) --(0,-5.25)|- (genhar2);
			\path (genhar2) -- (genh2);
			\path (genh2) -- (eigensolver2);
			\path (eigensolver2) -- (updaterho2);
			\path (updaterho2) -- (scfdec);
			\path (scfdec) -- node{No} (-4,-7)-- (-4,-1.5) -- (genhar);
			\path (scfdec) -- node{Yes}(4,-7)-- (adadec);
			\path (adadec) -- node{No}(8.75,-5)-| (adaks1);
			\path (adadec) -- (8.75,-5)-| (adaks2);
			\path (adadec) -- (8.75,-5)-| (adahar);
			\path (adaks1) |- (8.75,-2) |- (null1);
			\path (adaks2) |- (8.75,-2) |- (null1);
			\path (adahar) |- (8.75,-2) |- (null1);
			\path (adadec) -- node{Yes} (end);
		\end{scope}
		
	\end{tikzpicture}
	\caption{Flowchart of the multi-mesh adaptive algorithm for the KS equation. \label{fig:algorithm}}
\end{figure*}

Now that in the simulation, we have three meshes, for the wavefunctions, and the Hartree potential, respectively. The splitting algorithm for solving the KS equation is now ready and is illustrated in \Cref{fig:algorithm}.

\subsection{A soft-locking method}\label{subsec:soft-locking}
In the previous subsection, we addressed the splitting of eigenpairs, which leads to a generalized eigenvalue problem with eigenpairs expressed in two different approximate spaces established on two KS meshes. In this subsection, we will discuss how this problem is solved. The Kohn-Sham (KS) equation is discretized on two separate KS meshes, resulting in the following two subproblems:
\begin{subequations}
    \begin{align}
		A_1X_1 = \varepsilon M_1X_1,&A_1,M_1\in \mathbb{R}^{n_1\times n_1}, X_1 \in \mathbb{R}^{n_1\times p_1},\label{eq:mul-ksevp1}\\
		A_2X_2 = \varepsilon M_2X_2, &A_2,M_2\in \mathbb{R}^{n_2\times n_2}, X_2 \in \mathbb{R}^{n_2\times p}.\label{eq:mul-ksevp2}
  \end{align}
\end{subequations}
Here $n_1,n_2$ represents the number of Dofs for the respective KS meshes. The use of a multi-mesh adaptive technique reduces the size of the discretized eigenvalue problem compared to solving for all eigenvalues on a single mesh. Furthermore, the number of required eigenpairs on the first mesh is also smaller, enhancing computational efficiency. Specifically, for the first subproblem \eqref{eq:mul-ksevp1}, the goal is to compute the smallest $p_1$ eigenvalues and their associated eigenvectors. 

Unlike \Cref{eq:mul-ksevp1}, which can be solved directly using any iterative eigensolver such as the LOBPCG method \cite{knyazev2001toward, knyazev2007block, duersch2018robust} used in this work, the second subproblem \eqref{eq:mul-ksevp2} requires more careful handling. To compute the eigenpairs from the $(p+1)$-th to the $p$-th, the smallest $p_1$ eigenpairs are also needed, which is why $X_2$ contains $p$ eigenvectors. The first KS space is designed to accurately capture the first group of eigenvectors, while the second KS space targets the second group. Therefore, solving the second discretized eigenvalue problem requires accurate information from the first KS space's eigenpairs and an algorithm that focuses on finding the target eigenpairs defined in the second space, while also maintaining the orthogonality of the eigenvectors. 

We are motivated by the  \emph{soft-locking} strategy proposed in \cite{knyazev2007block} to design an effective approach for achieving the above goal. In LOBPCG, soft locking is used to lock converged eigenpairs during iterations, reducing computational cost. In our approach, we extend this concept to handle the first group of eigenpairs from the first KS space when solving the second discretized problem. By applying this strategy, we can focus computational efforts on the targeted eigenpairs in the second KS space. In the following we first briefly review the LOBPCG method \cite{knyazev2001toward} for solving the generalized eigenvalue problem, then the presented method follows.

\paragraph{LOBPCG}
Assume the problem reads as $A X =  BX\Lambda$. Here we denote the eigenvalues arranged in increasing order by $\lambda_1 \le \lambda_2 \le \cdots \le \lambda_n$. Their corresponding eigenvectors are denoted by $x_1,x_2,\cdots,x_n$. Then the first $k \leq n$ eigenvectors and eigenvalues are given by $X=\left[x_1, x_2, \ldots, x_k\right]$ and $\Lambda=$ $\operatorname{diag}\left\{\lambda_1, \lambda_2, \ldots, \lambda_k\right\}$, respectively. The eigenvectors are approximated by using the updating formula 
$$
X^{(i+1)}=X^{(i)} C_1^{(i+1)}+X_{\perp}^{(i)} C_2^{(i+1)},
$$
in the $i$-th iteration, where $X_{\perp}^{(i)}=\left[W^{(i)}, P^{(i)}\right]$.  The block $W^{(i)}$ is the preconditioned residual
\begin{equation}\label{eq:precon}
    W^{(i)}=T^{-1}\left(A X^{(i)}-B X^{(i)} \Lambda^{(i)}\right)
\end{equation}
with $\Lambda^{(i)}=X^{(i) T} A X^{(i)}$, where $T$ is any preconditioner. In this work, we use the preconditioner presented in \cite{bao2012h}. It has the form of $T = \frac{1}{2}L-\lambda B$, where $\frac{1}{2}L$ is the discretized kinetic operator in Hamiltonian, and $\lambda$ is an approximated eigenvalue. For each eigenpair we will construct a preconditioner to accelerate the calculation, as a  result, the preconditioners are designed as 
\begin{equation*}
T_{l}^{(i)}=\left\{\begin{aligned}
&\frac{1}{2} L-\lambda_l^{(i)} B & \text { if } \lambda_{l}^{(i)} <0, \\
&I & \text { otherwise, }
\end{aligned} \quad \text { for } l=1, \ldots, p\right.
\end{equation*}
In the practical simulations, the precondition process \Cref{eq:precon} involves the solution of a linear system, which can be implemented by using the algebraic multigrid (AMG) method. Specifically, there is no need to accurately solve this linear system, and few AMG iteration steps are performed in this process.

The block $P^{(i)}$ represents the $i$-th update direction, defined as
$$
P^{(i+1)}=X_{\perp}^{(i)} C_2^{(i+1)},
$$
with $P^{(1)}$ being an empty block, i.e., $X_{\perp}^{(1)}=W^{(1)}$. The coefficient matrices $C_1^{(i+1)}$ and $C_2^{(i+1)}$ are computed at each iteration of LOBPCG by the \textit{Rayleigh-Ritz} procedure, which involves solving a small generalized eigenvalue problem within the subspace $\mathcal{S}^{(i)}$ spanned by $X^{(i)}, W^{(i)}$, and $P^{(i)}$. Specifically,
\begin{equation}\label{eq:rrevp}
\left(S^{(i) T} A S^{(i)}\right) C^{(i+1)}=\left(S^{(i) T} B S^{(i)}\right) C^{(i+1)} \Lambda^{(i+1)},
\end{equation}
where $S^{(i)}$ is a matrix whose columns form a basis for $\mathcal{S}^{(i)}$, and is constructed as $S^{(i)}=\left[X^{(i)}, X_{\perp}^{(i)}\right]$. The corresponding matrix $C^{(i+1)}$ is given by
$$
C^{(i+1)}=\left[\begin{array}{ll}
	C_1^{(i+1)} & C_{1 \perp}^{(i+1)} \\
	C_2^{(i+1)} & C_{2 \perp}^{(i+1)}
\end{array}\right] .
$$
The leading $k$ columns of $C^{(i+1)}$ are
$$
C_x^{(i+1)}=\left[\begin{array}{l}
	C_1^{(i+1)} \\
	C_2^{(i+1)}
\end{array}\right],
$$
which are used to update $X^{(i+1)}$, while the remaining columns form the orthogonal complement within the search subspace.
	
The diagonal matrix $\Lambda^{(i+1)}$ in \eqref{eq:rrevp} contains approximations to the desired eigenvalues. If $k$ smallest eigenpairs are sought, the eigenvalues are sorted in ascending order.  

The soft-locking strategy in LOBPCG aims to deflate converged eigenvectors by keeping them in $X^{(i)}$ but exclude corresponding columns from $W^{(i)}$ and $P^{(i)}$. In order to produce reliable results, approximate eigenpairs (Ritz pairs) should be locked in order, i.e., the $(j+1)$-th Ritz pairs cannot be locked if the $j$-th Ritz pairs do not satisfy the convergence criterion.

\paragraph{Soft-locking in splitting method}
Motivated by the the soft-locking procedure, we propose an approach to solve the second subproblem \eqref{eq:mul-ksevp2} generated by the presented splitting method. The idea is similar: lock the converged eigenpairs and then iterated to find the unconverged ones. After solving the first discretized eigenvalue problem, we obtain the converged $p_1$ eigenpairs $\{(\varepsilon_l,X_l)\},l=1,\dots,p_1$. Then we adopt these eigenpairs to solve the second discretized eigenvalue problem in \Cref{eq:mul-ksevp2}. Two numerical challenges arise: first, expressing the first group of eigenvectors in the second KS space, given that the eigenvectors in these two spaces differ in size; second, preserving the orthogonality of the eigenvector across different spaces, which is nontrivial as orthogonality is difficult to define when they belongs to different discretized spaces.

The first issue is not particularly problematic, as locked eigenvectors evolve even in the solution of a single discretized generalized eigenvalue problem. We can simply interpolate the first group of eigenvectors from the first KS space to the second KS space and then iterate the LOBPCG method on the second KS space. The key point is that the first group of eigenvectors serves to accelerate the solution of the second discretized eigenvalue problem, and the final $p_1$ eigenpairs from the second KS space are not utilized. 

To address the second issue—ensuring that the eigenvectors (wavefunctions) are orthonormal despite originating from different finite element spaces—we apply a post-processing technique to enforce orthogonality. This process must be carried out in a unified finite element space. Specifically, we first construct a larger finite element space $V_h$, which encompasses all the wavefunctions by building it over the union of the meshes used for each wavefunction. The wavefunctions are then interpolated onto this unified space $V_h$. To mitigate interpolation errors, instead of directly re-orthogonalizing the wavefunctions, we solve a small generalized eigenvalue problem within the subspace spanned by the interpolated wavefunctions. The resulting eigenvectors from this small eigenvalue problem are the orthogonalized wavefunctions. 

In the practical simulation, the post-processing is performed only at the conclusion of the simulation. The construction of the finite element space $V_h$ can be efficiently handled using hierarchical grid transformation (HGT), and only a small eigenvalue problem needs to be solved, making this post-processing technique computationally efficient.

\subsection{Issues in the self-consistent field iteration}
We introduce the simple density mixing scheme related to the self-consistent field (SCF) iteration here. The traditional way to use the simple mixing is illustrated below. Assume the parameter $\alpha$ is the mixing factor, and the electron density at the $i$-th iteration is denoted by $\rho^{(i)}$. The electron density  at the $(i+1)$-th iteration is updated by
\begin{equation}\label{eq:mix}
	\rho^{(i+1)} = \alpha \rho^{(i)} + (1-\alpha)\rho^{new},
\end{equation}
where $\rho^{new}$ is the electron density obtained from the current SCF iteration. The mixing factor $\alpha$ is usually set to a small value, such as $0.618$ in this work. 

In the splitting algorithm, the electron density needs to be updated on each mesh, complicating the process slightly. Assume that in the $i$-th iteration, we have solved the first eigenvalue problem, and the wavefunctions are updated, denoted as $\{\psi^{(i+1)}_l\},~l=1,\dots,p_1$.  Before updating the electron density, these obtained wavefunctions are interpolated to other finite element spaces, denoted as 
\begin{align*}
	\text{On mesh $\mathcal{T}_2$:} &\{\psi^{(i+1)}_{l,2}\},~l=1,\dots,p_1,\\
	\text{On mesh $\mathcal{T}_H$:} &\{\psi^{(i+1)}_{l,H}\},~l=1,\dots,p_1.\\
\end{align*}
Then we can generate the new electron densities as 
\begin{equation}\label{eq:rho-new-ks-1}
	\begin{aligned}
		\rho^{new}_1 &= \sum_{l=1}^{p_1}f_l|\psi_l^{(i+1)}|^2 + \sum_{l=p_1+1}^{p}f_l|\psi_{l,1}^{(i)}|^2,\\
		\rho^{new}_2 &= \sum_{l=1}^{p_1}f_l|\psi_{l,2}^{(i+1)}|^2 + \sum_{l=p_1+1}^{p}f_l|\psi_l^{(i)}|^2,\\
		\rho^{new}_H &= \sum_{l=1}^{p_1}f_l|\psi_{l,H}^{(i+1)}|^2 + \sum_{l=p_1+1}^{p}f_l|\psi_{l,H}^{(i)}|^2.
	\end{aligned}
\end{equation}
Now we can perform the first mixing step for electron densities as 
\begin{equation}\label{eq:mix-ks-1}
	\begin{aligned}
		\tilde{\rho}^{i+1}_1 &= \alpha \rho_1^{(i)}+(1-\alpha)\rho^{new}_1,\\
		\tilde{\rho}^{i+1}_2 &= \alpha \rho_2^{(i)}+(1-\alpha)\rho^{new}_2,\\
		\tilde{\rho}^{i+1}_H &= \alpha \rho_H^{(i)}+(1-\alpha)\rho^{new}_h.
	\end{aligned}
\end{equation}
Next we update the Hartree potential and solve the second eigenvalue problem. After this, we can update the electron densities on the three meshes as follows
\begin{equation}\label{eq:rho-new-ks-2}
	\begin{aligned}
		\rho^{new}_1 &= \sum_{l=1}^{p_1}f_l|\psi_l^{(i+1)}|^2 + \sum_{l=p_1+1}^{p}f_l|\psi_{l,1}^{(i+1)}|^2,\\
		\rho^{new}_2 &= \sum_{l=1}^{p_1}f_l|\psi_{l,2}^{(i+1)}|^2 + \sum_{l=p_1+1}^{p}f_l|\psi_l^{(i+1)}|^2,\\
		\rho^{new}_H &= \sum_{l=1}^{p_1}f_l|\psi_{l,H}^{(i+1)}|^2 + \sum_{l=p_1+1}^{p}f_l|\psi_{l,H}^{(i+1)}|^2.
	\end{aligned}
\end{equation}
Finally we perform the second mixing step for electron densities as
\begin{equation}\label{eq:mix-ks-2}
	\begin{aligned}
		\rho^{i+1}_1 &= \alpha \tilde{\rho}_1^{(i+1)}+(1-\alpha)\rho^{new}_1,\\
		\rho^{i+1}_2 &= \alpha \tilde{\rho}_2^{(i+1)}+(1-\alpha)\rho^{new}_2,\\
		\rho^{i+1}_H &= \alpha \tilde{\rho}_H^{(i+1)}+(1-\alpha)\rho^{new}_H.
	\end{aligned}
\end{equation}

Using this mixing scheme, the electron densities on different meshes are updated to make the SCF iteration converge. A corresponding stop criterion for the SCF iteration can be delivered by the value of the difference between adjacent electron densities via
\begin{equation*}
\Delta \rho = \left\|\tilde{\rho}_1^{(i+1)}-\rho_1^{(1)}\right\|_2
+ \left\|\rho_2^{(i+1)}-\tilde{\rho}_2^{(i+1)}\right\|_2.
\end{equation*}
This formulation ensures the difference reflects the sum of the updates for each group of wavefunctions. The SCF iteration stops when $\Delta \rho$ falls below a predefined tolerance.

\section{Numerical examples}
In this section, we examine the convergence and efficiency of the splitting method through a series of numerical examples. In addition to the presented splitting strategy which is based on the splitting of core  and valence orbitals, we also explore two additional splitting strategies. All the simulations are performed on a workstation ``Moss" with two AMD EPYC 7713 64-Core Processors (at 2.0GHz$\times$64, 512M cache) and 900GB of RAM, and the total number of cores is 128. The software is the C++ library \texttt{AFEABIC} \cite{bao2012h, bao2015real, bao2016towards} under Ubuntu 20.04.

To measure the effectiveness of the splitting strategy, we propose a \textit{splitting factor}, defined by 
\begin{equation}
    \label{eq:split-factor}
    sf := \frac{\#\mathcal{T}_{\mathrm{KS},merged}}{\sum_{k=1}^{n_{grp}}\#\mathcal{T}_{\mathrm{KS},k}}.
\end{equation}
Here $\#\mathcal{T}_{\mathrm{KS},k}$ denotes the number of mesh grids for the $k$-th group for the eigenpairs, while $\#\mathcal{T}_{\mathrm{KS},merged}$ represents the number of grids in the merged mesh. It is noted that this merged mesh is obtained from merging the mesh for KS orbitals, not including the mesh for Hartree mesh, since in this work we focus on the splitting of the KS orbitals. For the discussions on the Hartree mesh, please refer to \cite{kuang2024towards}. From this definition \Cref{eq:split-factor}, the splitting factor ($sf$) indicates how effectively the KS problem is split, with a value between $1/n_{grp}$ and 1. Specifically, the closer $sf$ is to $1$, the more distinct the behaviors between different groups, indicating better splitting performance.

\subsection{Illustration of the splitting method: LiH molecule}
We begin by examining the method for solving the ground state of the LiH molecule. LiH contains $N_\mathrm{ele}=4$ electrons, and the system requires the solution of two occupied Kohn--Sham orbitals ($p=2$). To implement the splitting method, we simply divide these two orbitals into two groups, i.e., $p_1=p_2=1$. Consequently, the splitting method facilitates the adaptive generation of three finite element spaces: one for the Hartree potential and two for the two occupied orbitals, respectively.

To validate the accuracy of our approach, we compare the total energy obtained using our method with those from the state-of-the-art DFT software \texttt{NWChem} \cite{valiev2010nwchem}. We further compare our results with the adaptive finite element method on a single mesh to demonstrate the efficiency of our method. The results are summarized in \Cref{tab:LiHTime}. Specifically, we set the error indicator tolerance to $tol_{ada}=4e^{-6}$ for both methods. The table shows that our method achieves chemical accuracy (i.e., an energy error per atom of less than 0.001 Hartree), demonstrating the effectiveness of the splitting method. Although the splitting method achieve a similar accuracy with the single mesh method, it is worthy to note that the size of the eigenvalue problem is reduced by two thirds. This reduction means that instead of solving a single eigenvalue problem with a matrix size of $6,029,109\times 6,029,109$, we solve two smaller eigenvalue problems with sizes of $2,132,602\times 2,132,602$ and $1,997,052\times 1,997,052$, respectively. From this point of view, the computational cost is reduced. 
\begin{table}[!htp]
	\caption{Comparison on \ch{LiH}. The referenced total energy is $E_\mathrm{ref}=-7.919582$ Hartree from \texttt{NWChem} using \emph{aug-cc-pv5z} basis set. \label{tab:LiHTime}}
	\centering
	\begin{tabular}{rrrr}
		\toprule
		&single mesh& split method & split method(+) \\ \midrule
		$E_{tot}$ &-7.917843 &-7.917805 &  -7.918257 \\
		$\Delta E/2$ &0.00087  & 0.00089 &0.00066 \\
		\#$\mathcal{T}_{\mathrm{KS},1}$ &
		6,029,109 &2,132,602 & 3,780,288 \\
		\#$\mathcal{T}_{\mathrm{KS},2}$ & -&1,997,052 &-\\
		\#$\mathcal{T}_{\mathrm{Har}}$ & - & 3,664,979 & 3,664,979\\
		\bottomrule  
	\end{tabular}	
\end{table}

In addition, the accuracy is improved after we do the post-processing procedure (\emph{split method(+)}), as indicated in the right column of \Cref{tab:LiHTime}. {The effectiveness of the splitting method can also be justified by the aforementioned splitting factor. In this example, $sf=0.915$ which is close to $1$, hence the splitting strategy works fine in this example.}

\begin{figure}[!htp]
	\centering
	\includegraphics[width=0.85\linewidth]{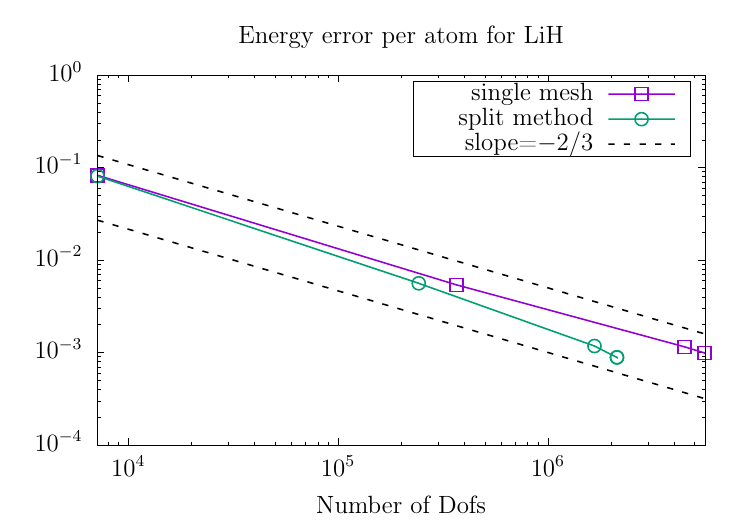}
	\includegraphics[width=0.85\linewidth]{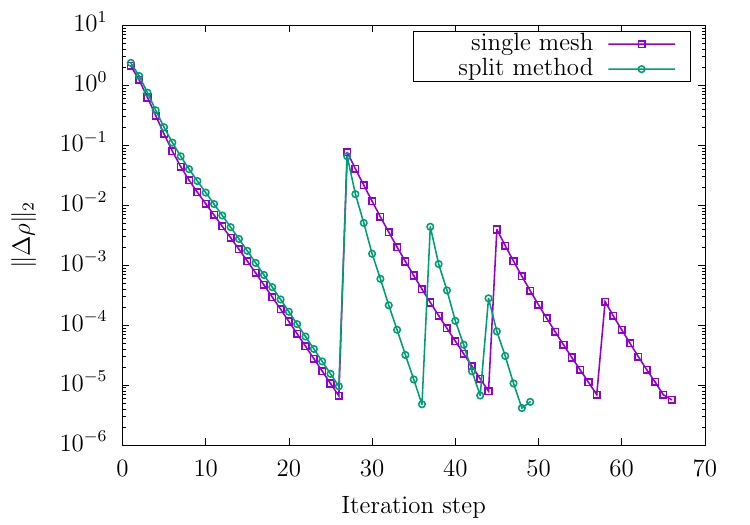}
	\caption{SCF iterations for LiH molecule. Top: convergence of energy for the single mesh adaptive method and the splitting method. Bottom: convergence of the electron density. \label{fig:lihscf}}
\end{figure}

We futher check the SCF convergence for the LiH molecule, as shown in \Cref{fig:lihscf}. The plot on the left illustrates the convergence of the total energy with respect to the number of mesh grids. From this figure, we observe systematic convergence of the total energy for both methods, with a rate consistent with the theoretical prediction, i.e., $\mathcal{O}(N^{-2/3})$ \cite{ciarlet2002finite}. Additionally, the splitting method demonstrates a slightly higher convergence rate compared to the single mesh method. The plot on the right depicts the convergence of the electron density. It is important to note that each increment in $\Delta \rho$ indicates a mesh adaptation. After the first mesh adaptation, the splitting method shows a more rapid decrease in $\Delta \rho$ compared to the single mesh method. This can be attributed to the fact that the density update is performed twice in each SCF iteration, as two eigenvalue problems are solved in each iteration. Consequently, the splitting method achieves faster convergence. These results highlight the higher accuracy and quicker convergence of the splitting method relative to the single mesh method under the same conditions.

\begin{figure*}[!htp]
	\centering
	\includegraphics[width=0.195\linewidth]{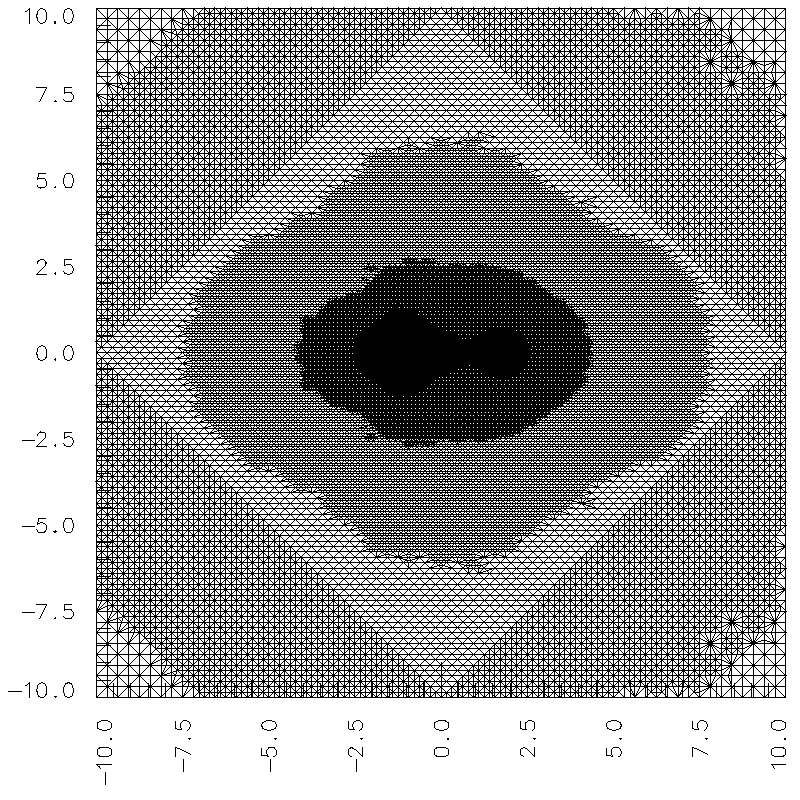}	
	\includegraphics[width=0.195\linewidth]{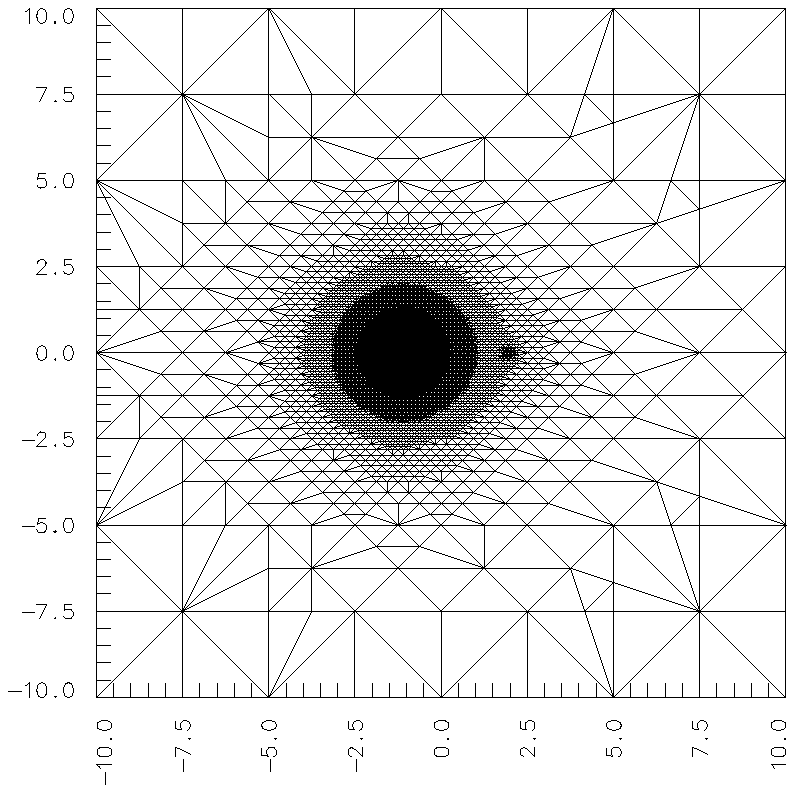}
	\includegraphics[width=0.195\linewidth]{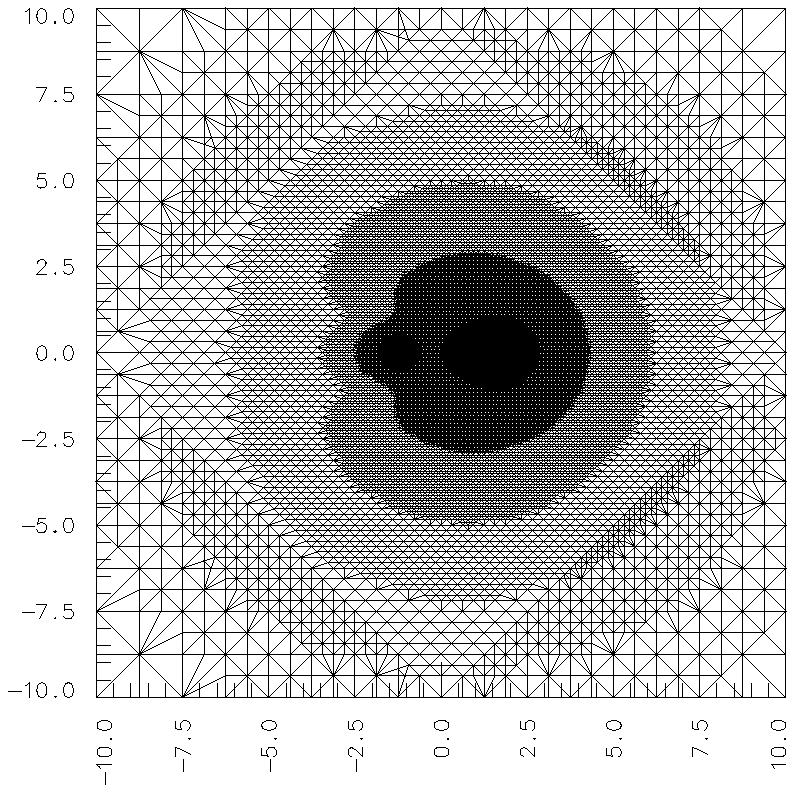}
	\includegraphics[width=0.195\linewidth]{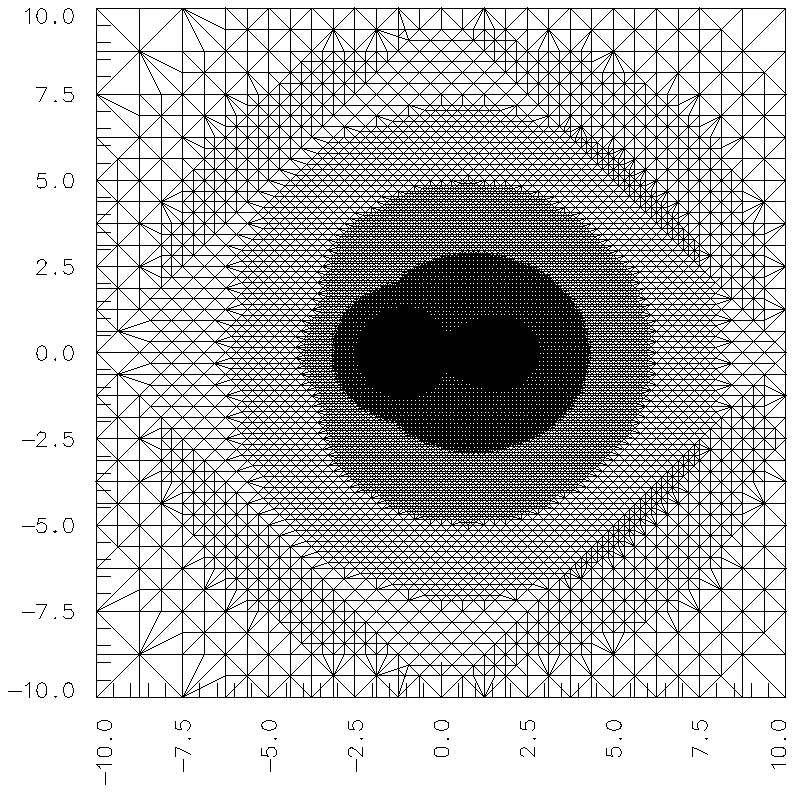} 
	\includegraphics[width=0.195\linewidth]{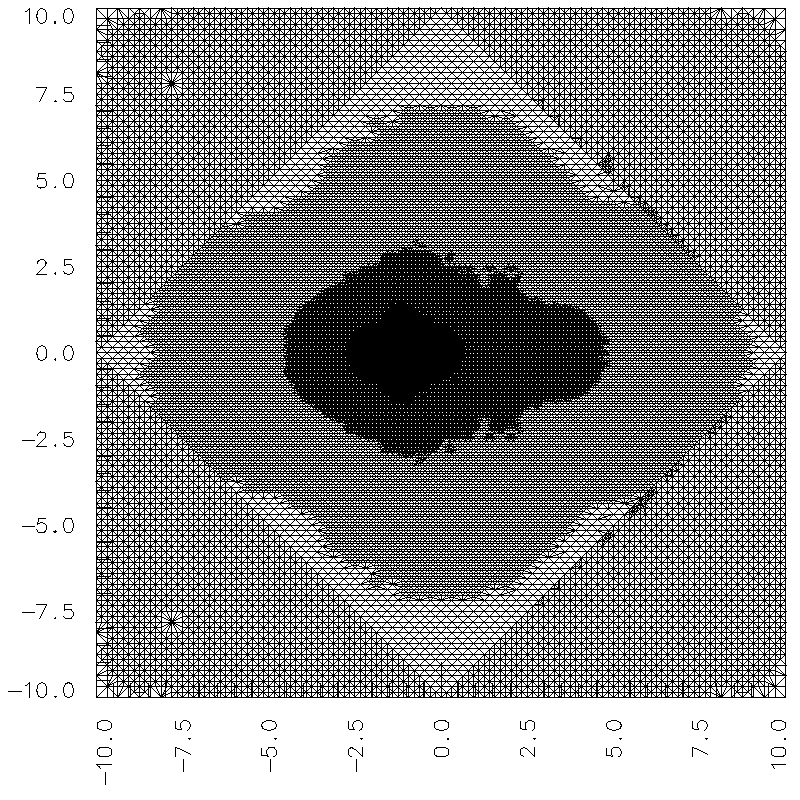}
	\includegraphics[width=0.195\linewidth]{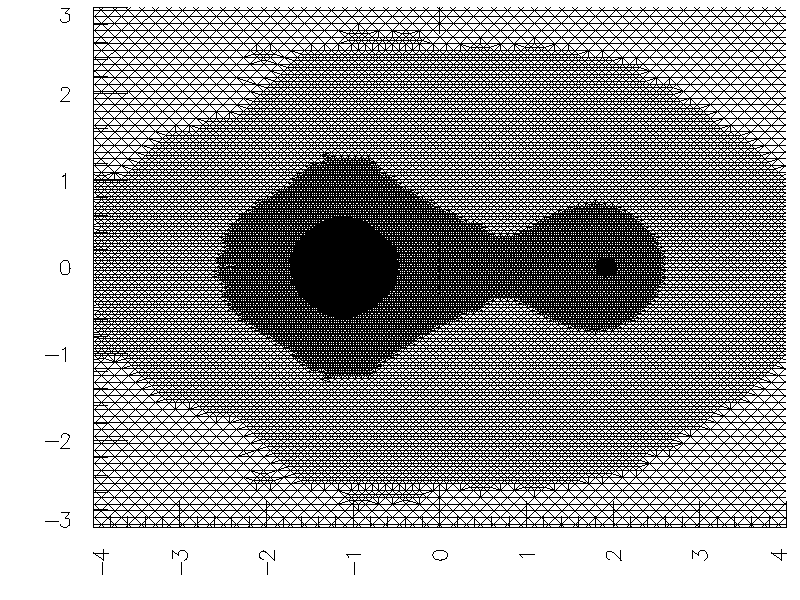}
	\includegraphics[width=0.195\linewidth]{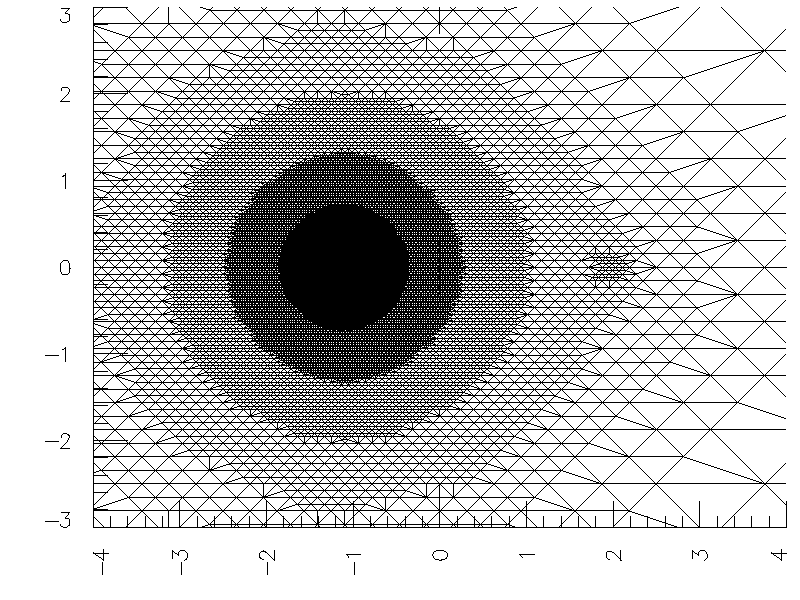}
	\includegraphics[width=0.195\linewidth]{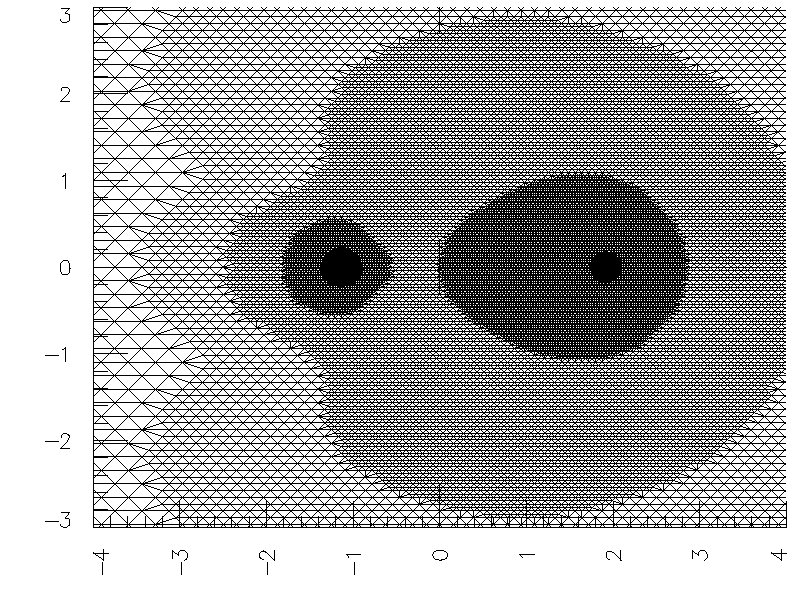}
	\includegraphics[width=0.195\linewidth]{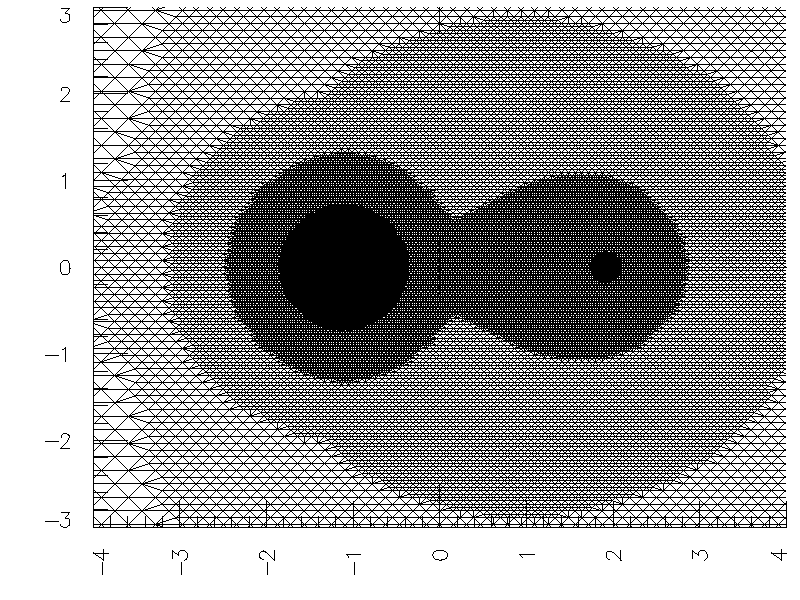}
	\includegraphics[width=0.195\linewidth]{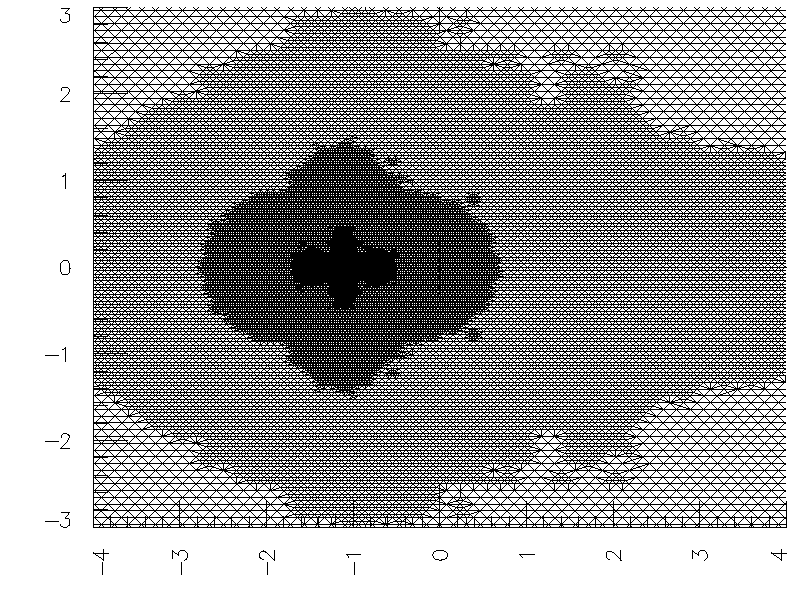}
	\includegraphics[width=0.195\linewidth]{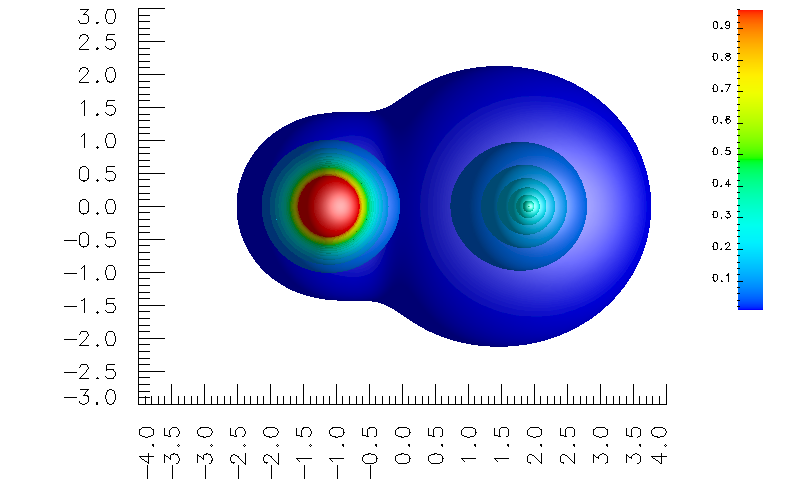}
	\includegraphics[width=0.195\linewidth]{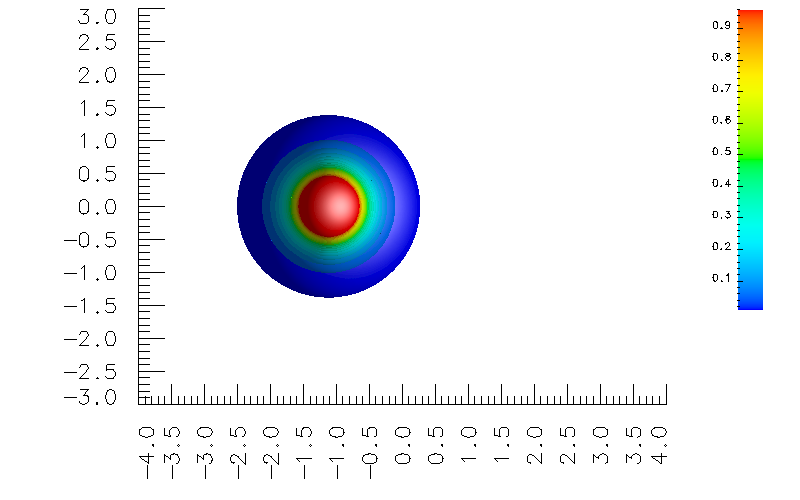}
	\includegraphics[width=0.195\linewidth]{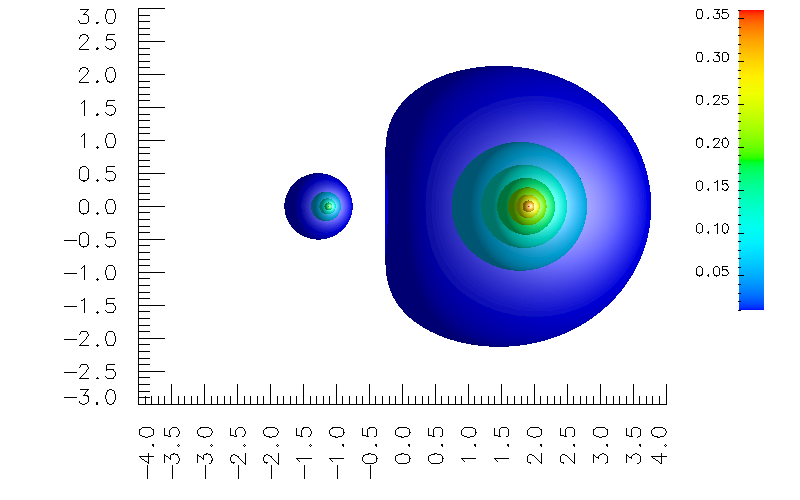}
	\includegraphics[width=0.195\linewidth]{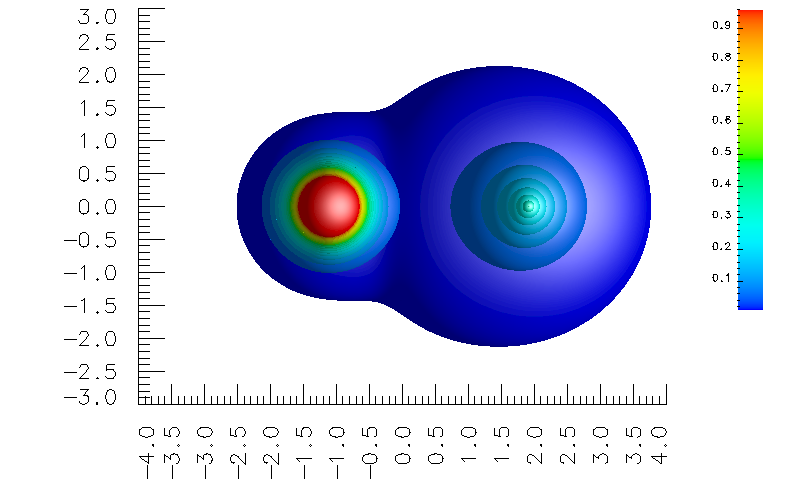}
	\includegraphics[width=0.195\linewidth]{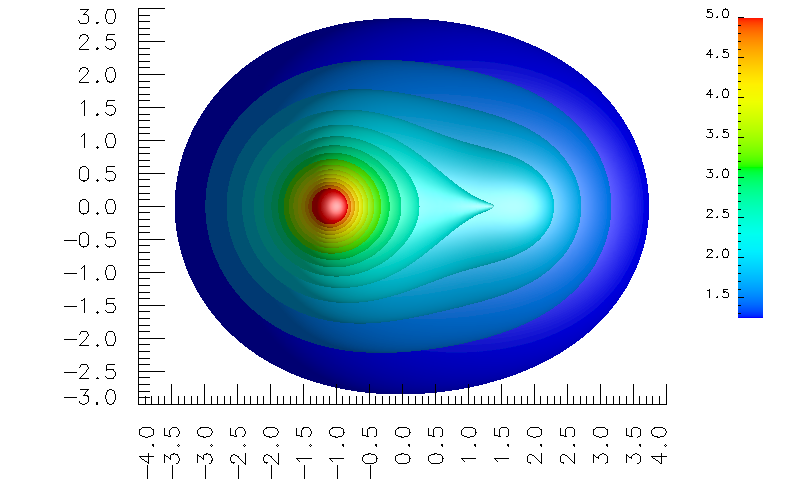}
	\caption{Left column is the results using single mesh method: sliced mesh on $X-Y$ plane (above), zooming-in for the sliced mesh (middle), and the electron density profile (below). Right four columns display the results using split-ks method: sliced meshes on $X-Y$ plane, zooming-in for the sliced mesh (middle), and associated electron density or Hartree potential profiles for  $\mathcal{T}_{\mathrm{KS},1}$, $\mathcal{T}_{\mathrm{KS},2}$, $\mathcal{T}_{\mathrm{KS},merged}$, and $\mathcal{T}_{\mathrm{Har}}$ from left to right. \label{fig:lih}}
\end{figure*}

For a more intuitive demonstration, we have plotted the meshes and contours for both the splitting method and the single method in \Cref{fig:lih}. The first column displays the results from the single method, while the remaining four columns show the results from the splitting method. This figure visually highlights the mechanics of the splitting method. In the single mesh method, the mesh shown in the top left of \Cref{fig:lih} is adapted to well capture both Kohn--Sham orbitals and the Hartree potential, resulting in a finer mesh near the nuclei. Specifically, the lithium atom is positioned at $(-1.11,0,0)$, and the hydrogen atom at $(1.92,0,0)$, leading to a denser grid in these regions. 

In contrast, the splitting method tailors each mesh to capture the variations of the corresponding variable. The second column of \Cref{fig:lih} shows that the first KS orbital which corresponds to the core electrons mainly varies around the lithium atom, resulting in a dense grid around lithium and a sparser one near hydrogen. The third column illustrates the second KS orbital which corresponds to the valence electrons, which shows significant variation in a large region around the hydrogen atom and a smaller region around lithium. After applying the post-processing procedure, the fourth column displays the merged mesh that effectively captures variations in all the KS orbitals. The fifth column of \Cref{fig:lih} shows the mesh grid distribution for the Hartree potential, which is denser near the boundaries and sparser near the nuclei compared to the KS meshes.

The grid distributions align with the shapes of the corresponding electron density and Hartree potential,, as shown in the last row of \Cref{fig:lih}. Notably, the splitting method simplifies the solution of the eigenvalue problem on the original dense mesh (left column) by splitting it into two problems on the smaller, less dense meshes (second and third columns). Additionally, in the splitting method, the Hartree potential is solved on a separate mesh (fifth column), whereas in the single mesh method, it is solved on the same mesh as the KS orbitals. Overall, these issues make the splitting method more efficient than the single mesh method in solving the eigenvalue problem and the linear system for the Hartree potential.

\subsection{Efficiency of the splitting method: \ch{BeH2} molecule}
\begin{figure}[!htp]
	\centering
	\includegraphics[width=0.323\linewidth]{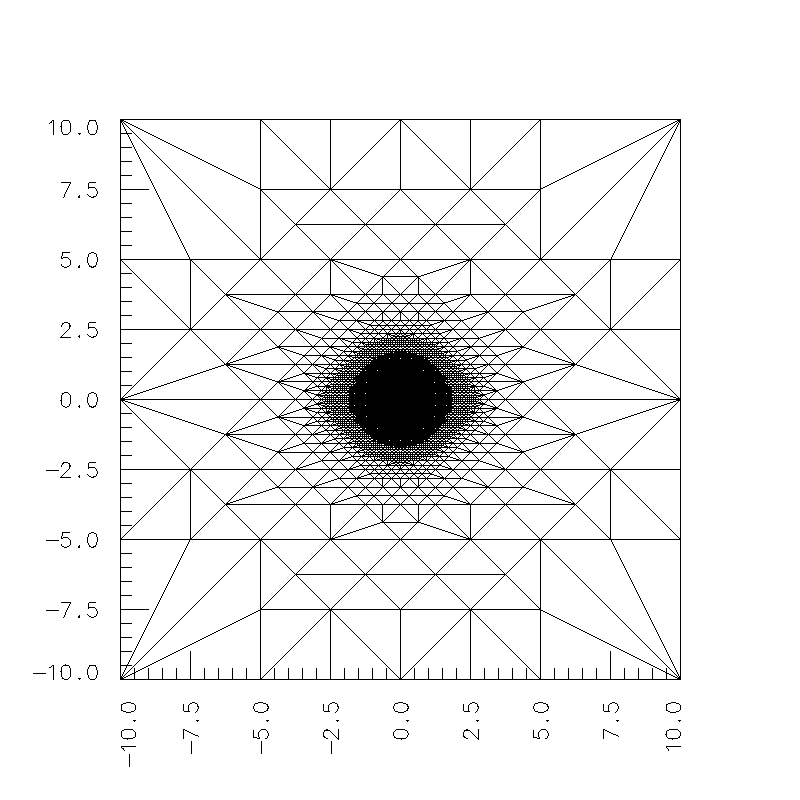}
	\includegraphics[width=0.323\linewidth]{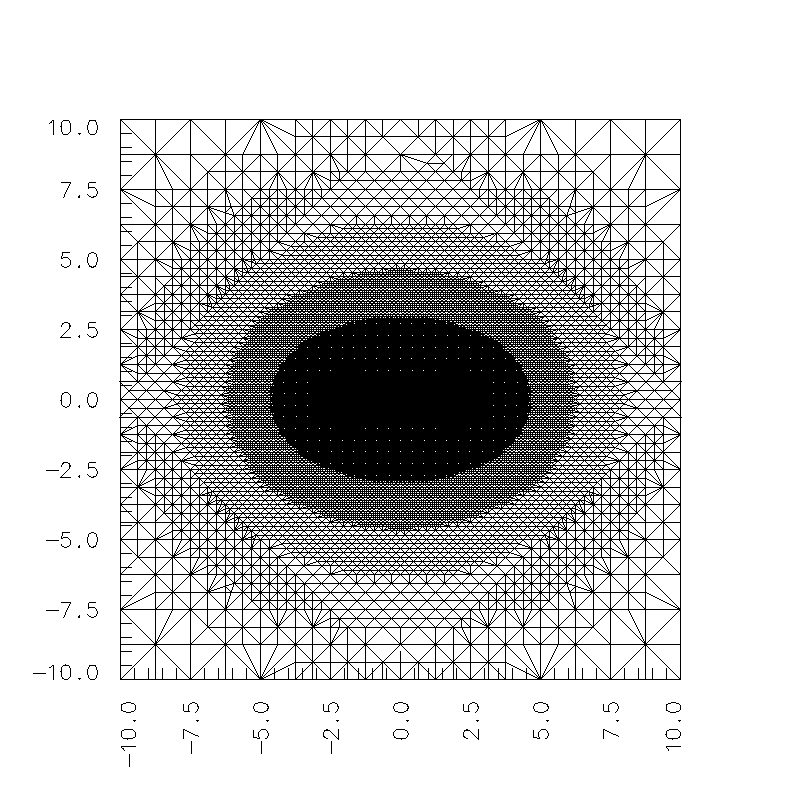}
	\includegraphics[width=0.323\linewidth]{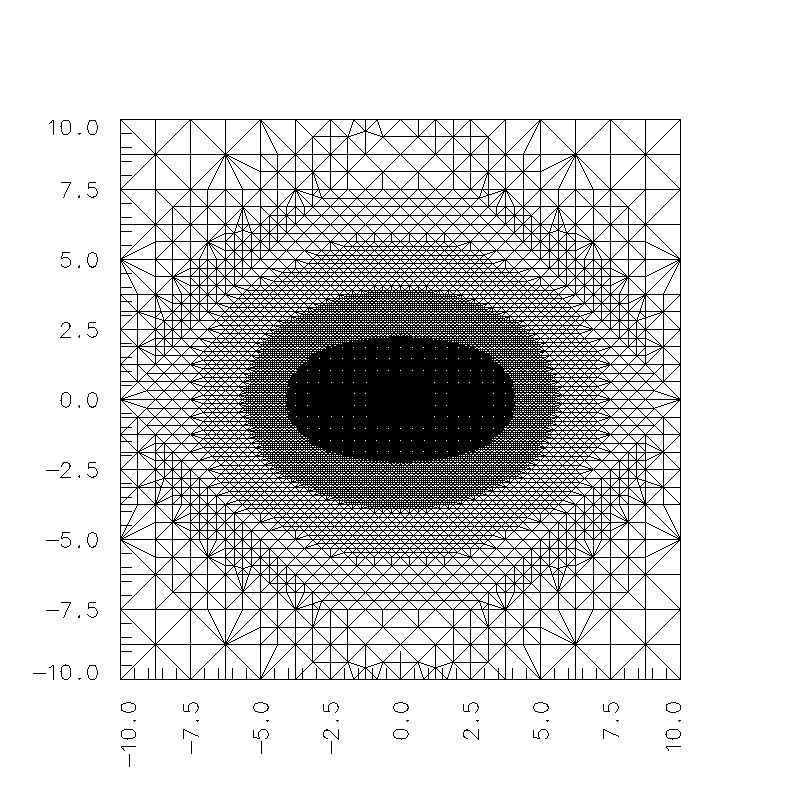} 
	\includegraphics[width=0.323\linewidth]{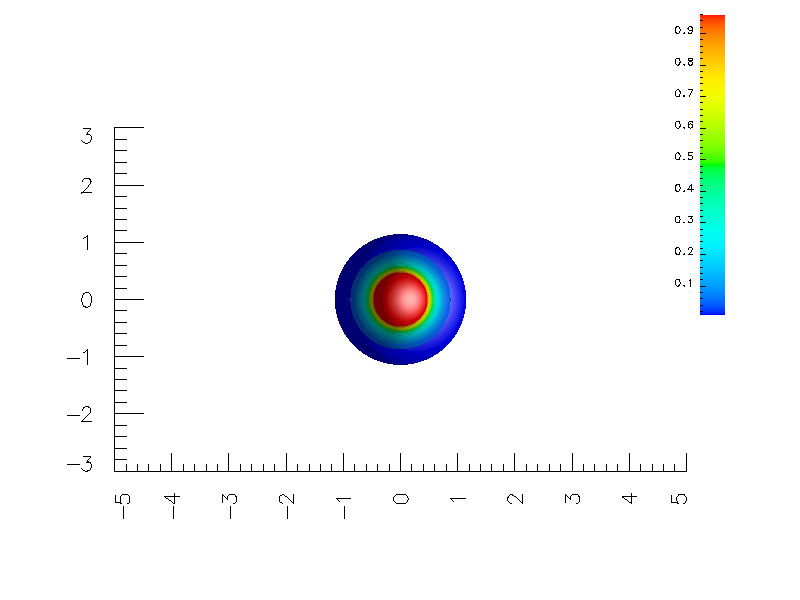}
	\includegraphics[width=0.323\linewidth]{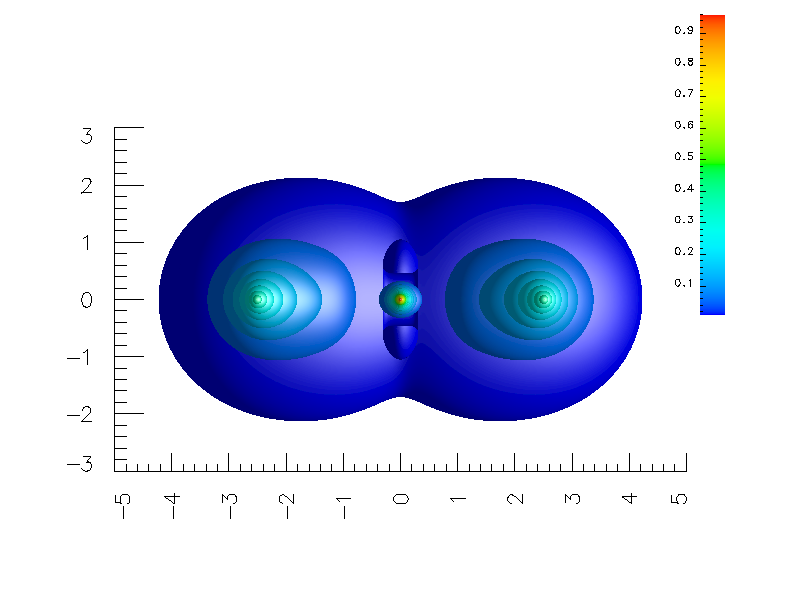}
	\includegraphics[width=0.323\linewidth]{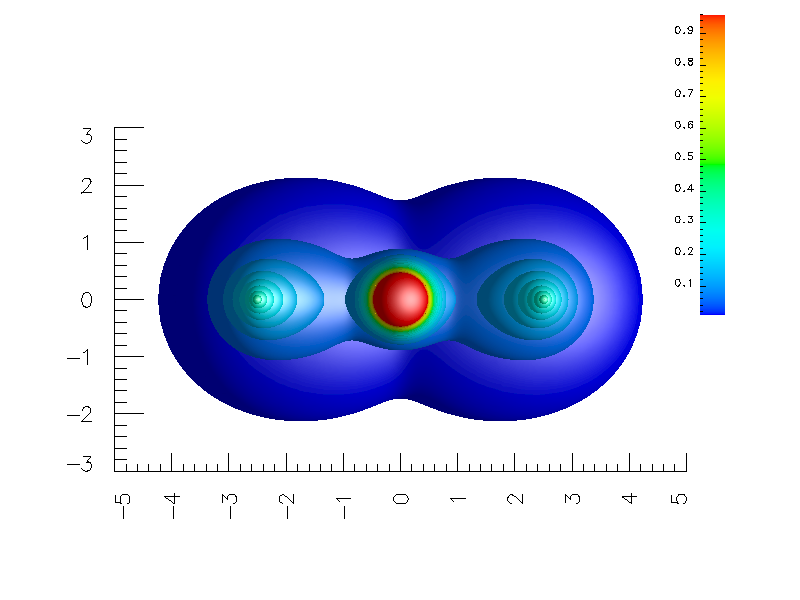}	
	\caption{Left: multimesh-har version, 4,852,528 DOFs on KS mesh, energy error $8.9e^{-4}$. Right two: multimesh-ks version, 2,324,465 DOFs on first KS mesh and 3,061,932 DOFs on second KS mesh,  energy error $1.1e^{-3}$. Bottom left: contour of electron density. Bottom middle: contour of electron density for orbitals in the first mesh. Bottom middle: contour of electron density for orbitals in the second mesh.   \label{fig:beh2}}
\end{figure}

To futher examine the efficiency and accuracy of the presented method, we then apply the splitting method to the \ch{BeH2} molecule. As indicated in \Cref{tab:split-mole}, the \ch{BeH2} molecule contains $N_\mathrm{ele}=6$ electrons, and the system requires the solution of $p=3$ occupied Kohn--Sham orbitals. We set $p_1=1, p_2=2$ for the splitting method based on the strategy which separate the orbitals into inner group and valence group. The meshes and contours for the \ch{BeH2} molecule are displayed in \Cref{fig:beh2}. Similar with the results for LiH, the mesh for different eigenpair group behaves different. Specifically, the first group corresponding to the core orbital has a denser mesh near the beryllium atom on a relatively small region, while the second group corresponding to the valence electrons has a denser mesh near all atoms on a relatively larger region. This demonstrates the effectiveness of the splitting.

\begin{figure}
	\centering
	\includegraphics[width=0.85\linewidth]{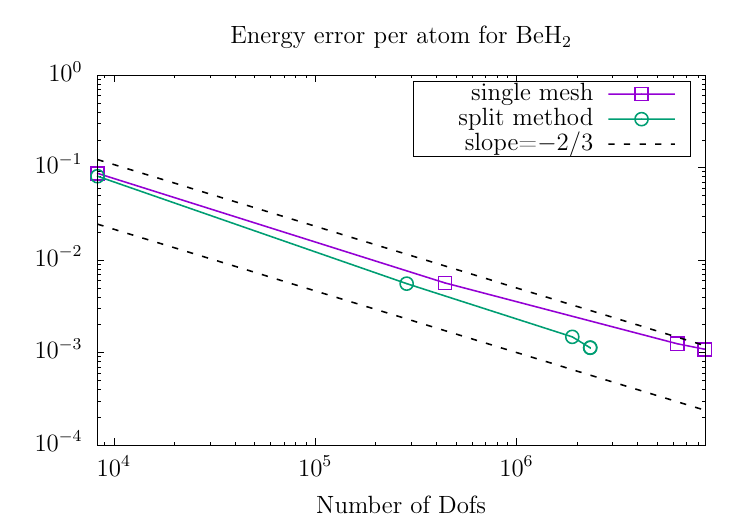}\\
	\includegraphics[width=0.85\linewidth]{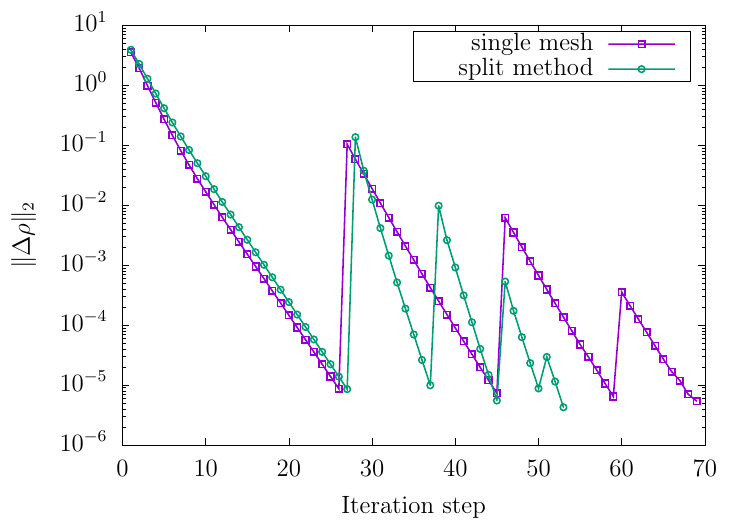}
	\caption{Error of total enery for molecule \ch{BeH2}.\label{fig:beh2conv}}
\end{figure}

The SCF convergence for the \ch{BeH2} molecule is shown in \Cref{fig:beh2conv}. Similar conclusions to the LiH molecule can be delivered  from this figure, which shows the generality of the presented method for different examples.

\begin{figure*}[!htp]
	\centering
	\includegraphics[width=\linewidth]{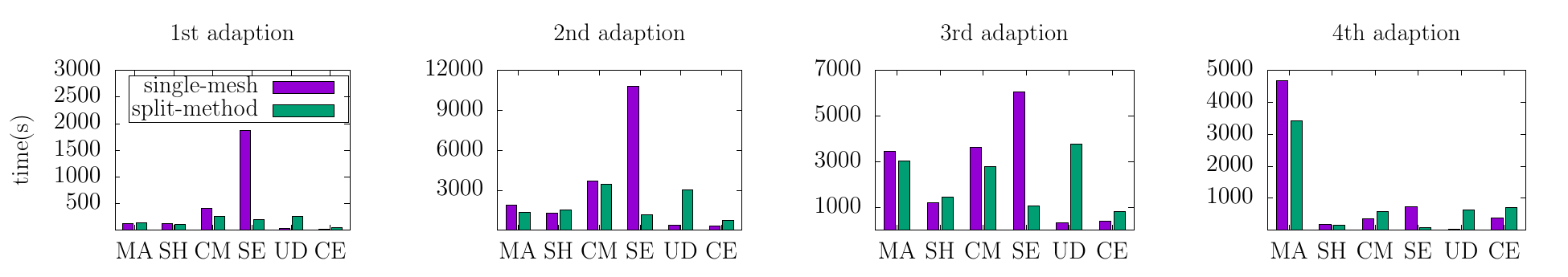}
	\caption{ CPU time for single mesh method and splitting method on six parts:  \textbf{MA} (mesh adaption), \textbf{SH} (solve Hartree potential), \textbf{CM} (construct matrix), \textbf{SE} (solve eigenvalue problem), \textbf{UD} (update electron density) and \textbf{CE} (calculate energy). Four mesh adaptations are needed to achieve chemical accuracy in both methods starting from the same initial setup. \label{fig:beh2-dtime}}
\end{figure*}

We now focus on the efficiency comparing to the single mesh method. The results are summarized in \Cref{tab:BeH2Time}. The table shows that the splitting method demonstrates a similar accuracy to the single mesh method, but with a reduced computational cost. From the table, more than $25\%$ CPU time is saved. It is worthy to mention that instead of solving a single eigenvalue problem with a matrix size of $7,338,638\times 7,338,638$, we solve two smaller eigenvalue problems with sizes of $2,631,201\times 2,631,201$ for the core orbital and $3,454,554\times 3,454,554$ for the valence orbitals, respectively.  The reduction in the size of the eigenvalue problem is significant, which leads to a reduction in the computational cost as shown in \Cref{fig:beh2-dtime}. After the post-processing process, the accuracy is also improved. {The splitting factor is $sf=0.837$ which also shows the effectiveness of the splitting strategy.}

In \Cref{fig:beh2-dtime}, we compare the CPU time in six parts according to the routine process: mesh adaption (MA), solution of Hartree (SH),  construction of matrices(CM), solution of the eigenvalue problems (SE), update electron density (UE), and calculation of total energy (CE).  The four images represent CPU time on four meshes during the mesh adaption process. It is observed that the MA part and CM part are reduced slightly, this is because two meshes and two eigenvalue problems are required in the splitting method. THen we check the CPU time in solving the eigenvalue problems (SE) and the linear system (SH part). From the detailed comparison in CPU time in \Cref{fig:beh2-dtime}, we find that the computational cost for the slution to the eigenvalue problem is significantly reduced, which is in accordance with the reduction of the matrix size. For SH part, it is noted that the cost is not reduced, this is because the generation of the right hand side of the linear system is involved in this part and it requires the information of the electron density, which implies that we need to communicates  from the Hartree mesh with  all the KS meshes to generate the electron density. Hence the cost of the linear system for the Hartree potential is not reduced. The CE part is also not reduced since the same reason as  the SH part. Finally, it is noted that the UD part is increased a lot, this is because the electron density is updated twice in each SCF iteration, and each time we need perform the interpolations of the wavefunctions for the soft lokcing procedure, while this requires a lot communications among meshes. Nevertheless, the total CPU time is reduced, which is mainly due to the reduction of the size of the eigenvalue problem.
\begin{table*}[!htp]
	\caption{Comparison on \ch{BeH2} with respect to results and serial computational time. The referenced value is  $E_\mathrm{ref} = -15.660638$ Hartree. \label{tab:BeH2Time}}
	\centering
	\begin{tabular}{ rrrrrrrrr}
		\toprule
		&$E_{tot}$ &$\Delta E/3$ &\#$\mathcal{T}_{\mathrm{KS},1}$ &\#$\mathcal{T}_{\mathrm{KS},2}$&\#$\mathcal{T}_{\mathrm{Har}}$& $t_{all}$& $t_{SCF}$ &$t_{MA}$\\
		\midrule
		single mesh&-15.657304& 0.00111&7,338,638 &- &-&42445.26 & 32331.87
		&10113.40\\\midrule
		split method&-15.657595&0.00101 & 2,631,201& 3,454,554 &4,741,565&30757.65&{22802.61}&{7955.04}\\
		split method(+) & -15.657947&0.00089&5,090,969&-&4,741,565& -& 355.76 & 3428.27\\
		\bottomrule  
	\end{tabular}	
\end{table*}

\subsection{Splitting based on eigenvalues: \ch{H2O} molecule}
In the above two examples, we adopt the splitting strategy which split the KS orbitals into inner core orbitals and the valence orbitals. The results have demonstrate the validity and efficiency for this splitting strategy. In this example, we would like to discuss another way of the splitting method in a \ch{H2O} example. The idea is straightforward. The regularities for the KS orbitals are expected to be different if they correspond to difference eigenvalue. 

\begin{figure}[!h]
	\centering
	\includegraphics[width=0.3\linewidth]{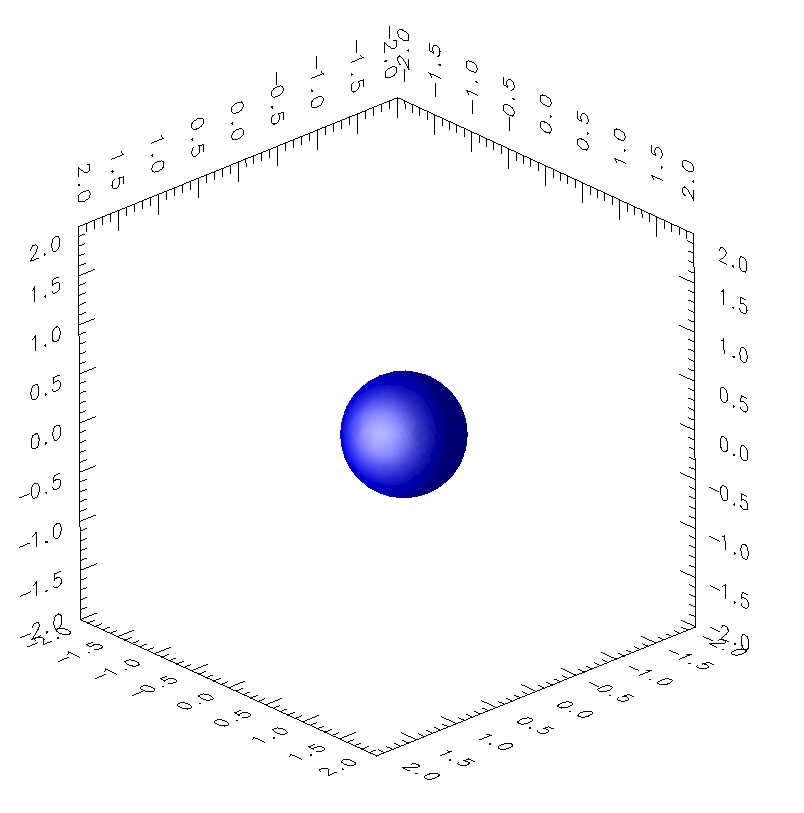}
	\includegraphics[width=0.3\linewidth]{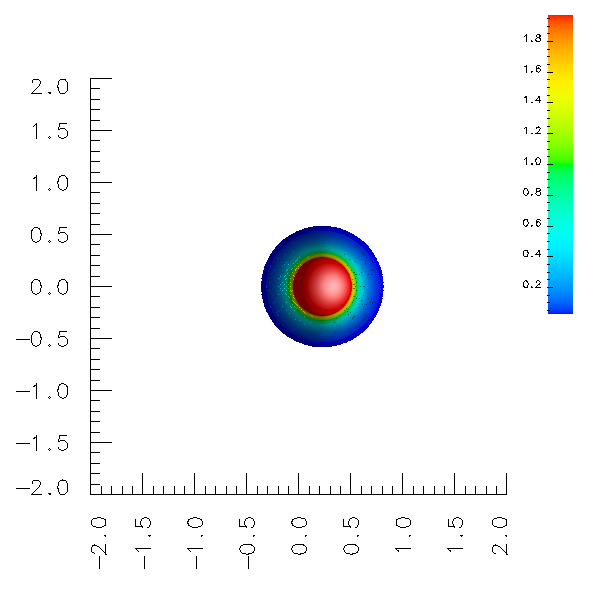}
	\includegraphics[width=0.3\linewidth]{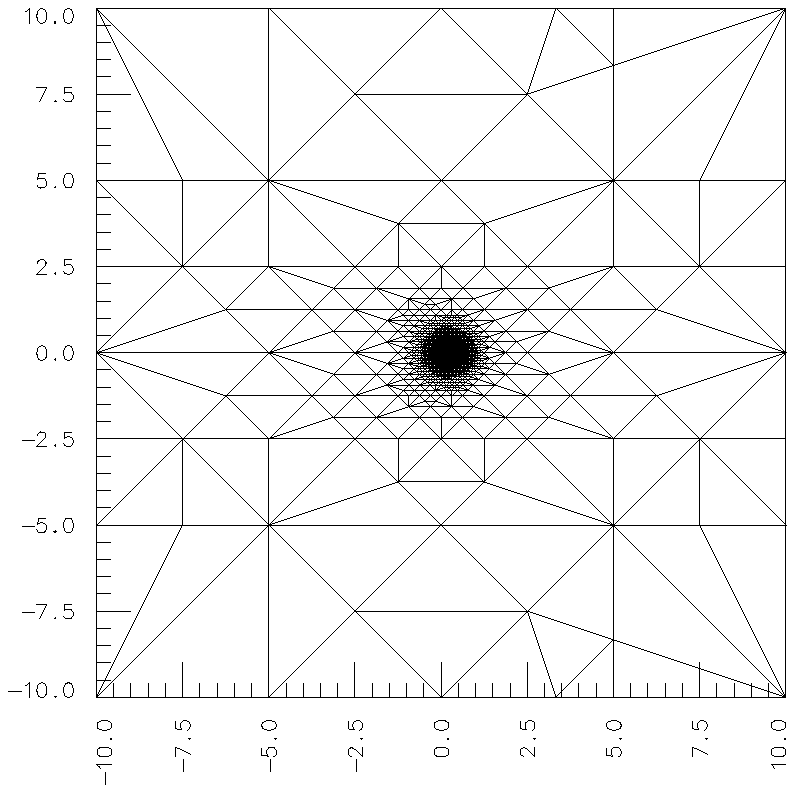}\\
	\includegraphics[width=0.3\linewidth]{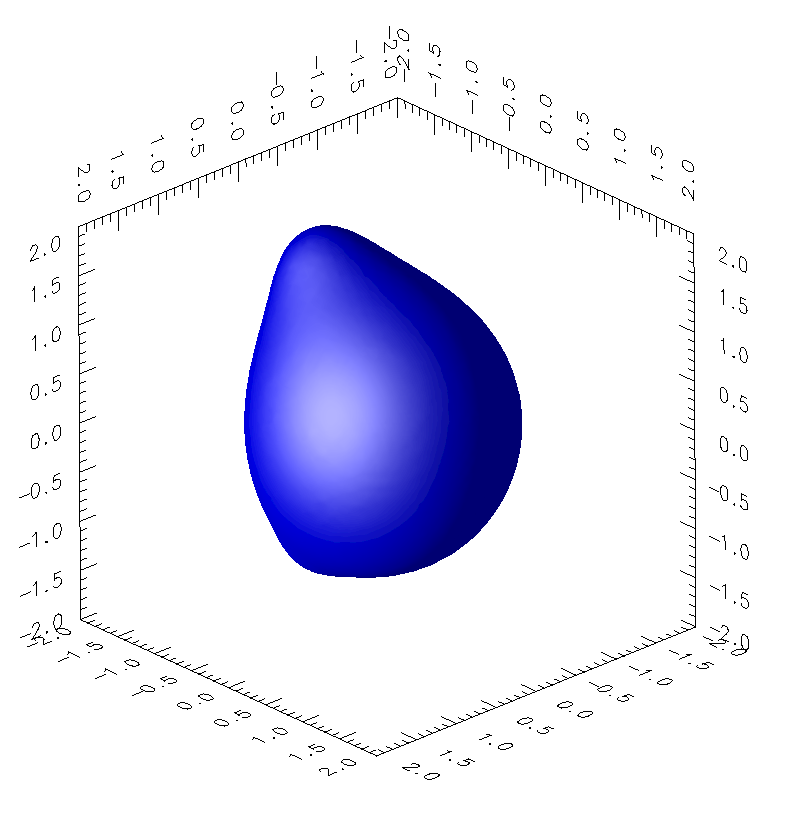}
	\includegraphics[width=0.3\linewidth]{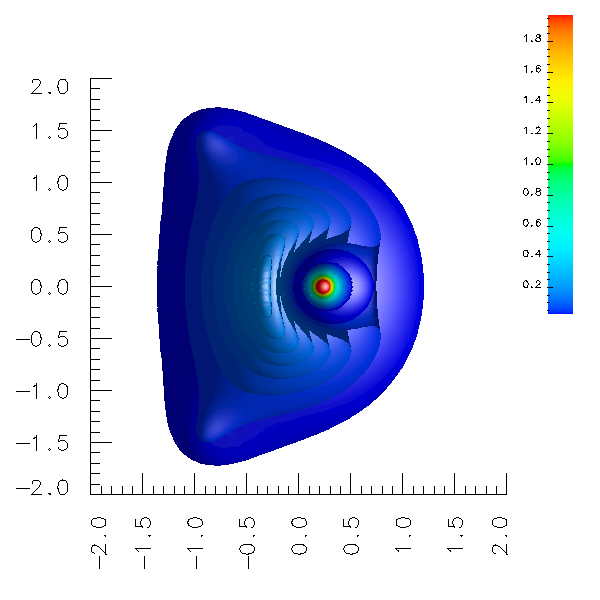}
	\includegraphics[width=0.3\linewidth]{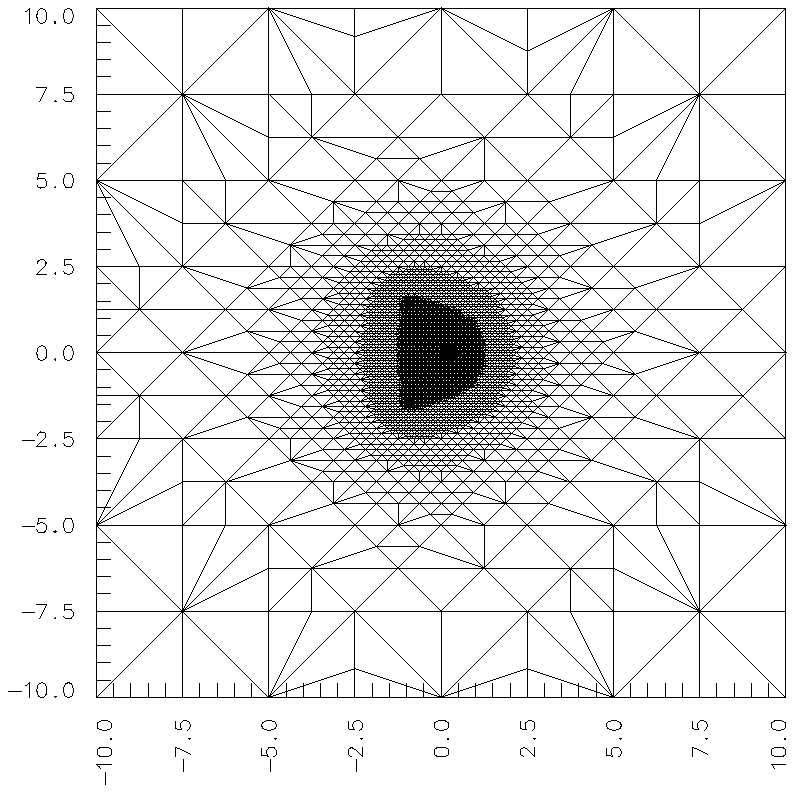}\\
	\includegraphics[width=0.3\linewidth]{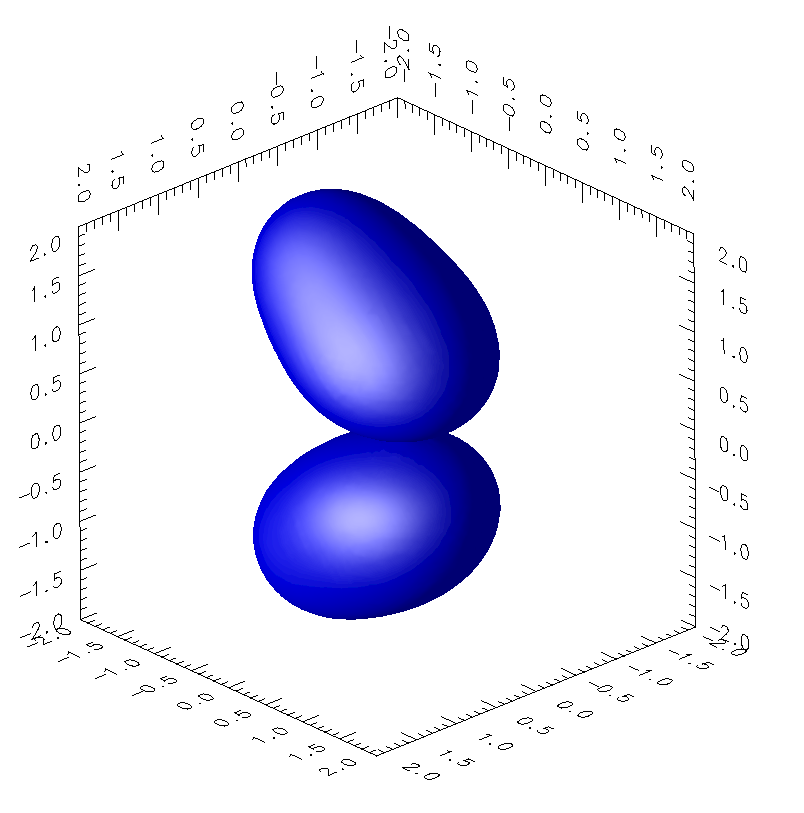}
	\includegraphics[width=0.3\linewidth]{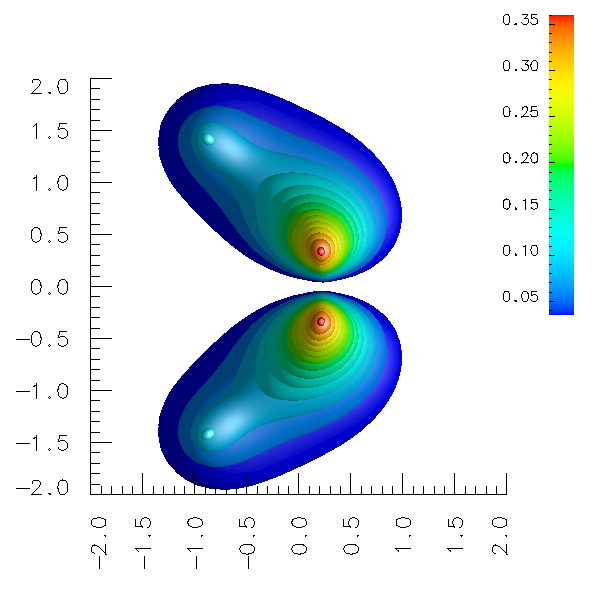}
	\includegraphics[width=0.3\linewidth]{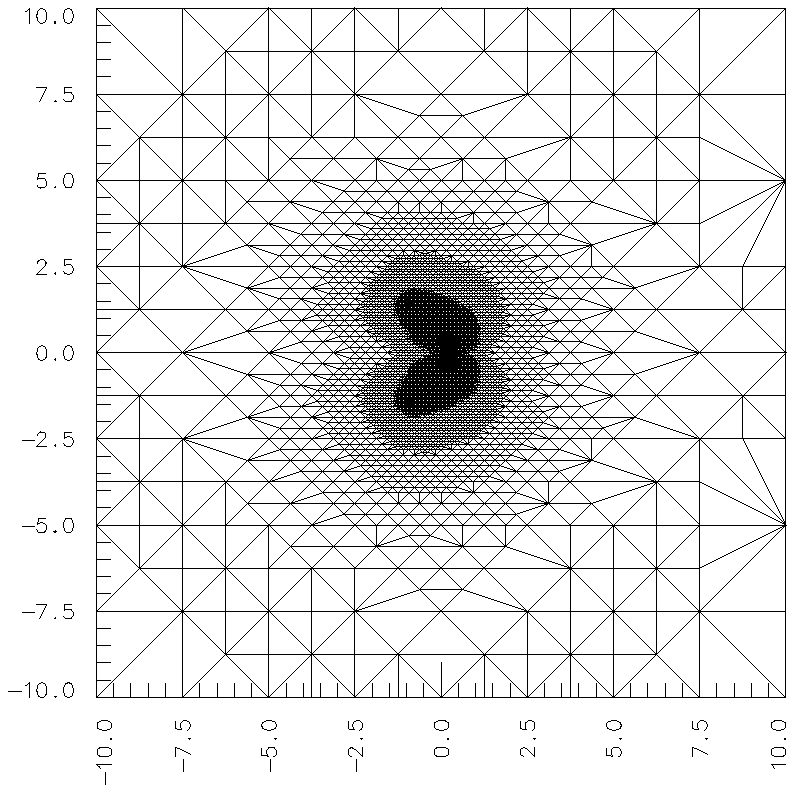}\\
	\includegraphics[width=0.3\linewidth]{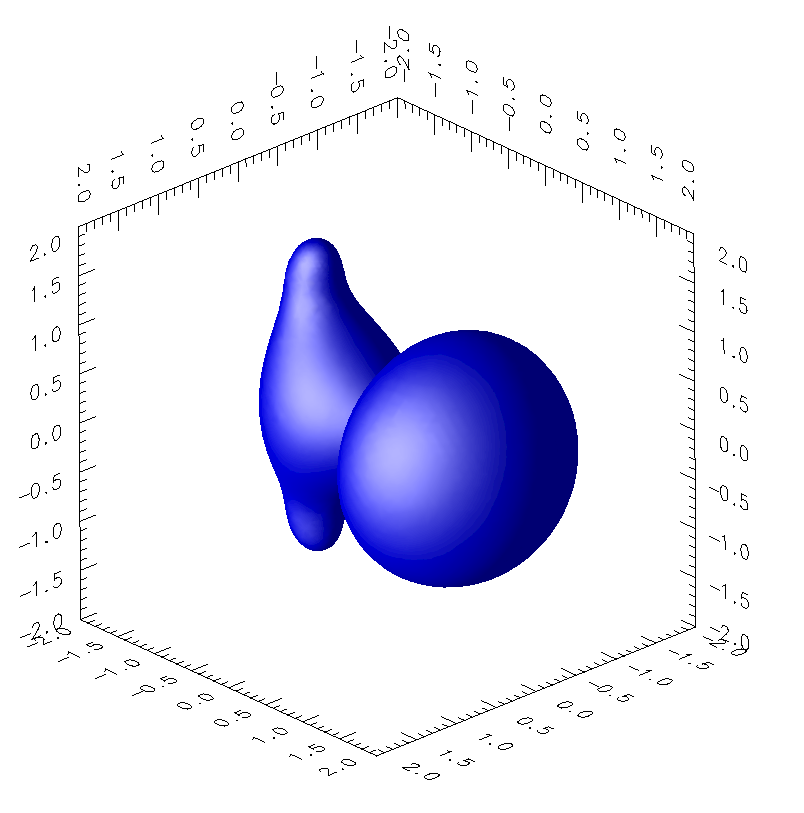}
	\includegraphics[width=0.3\linewidth]{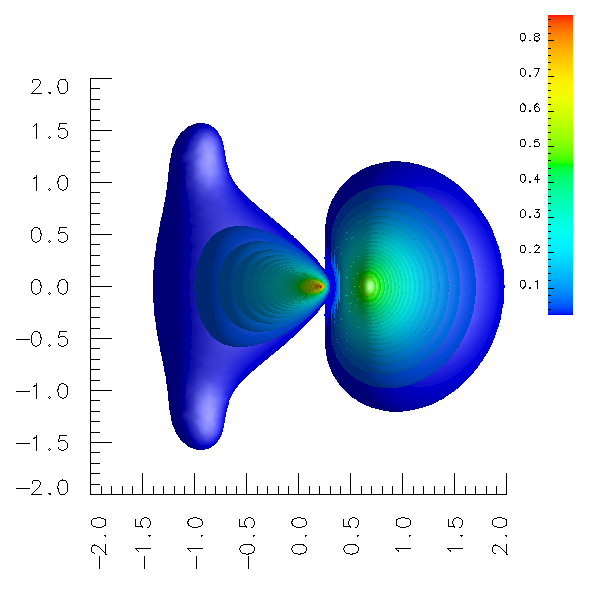}
	\includegraphics[width=0.3\linewidth]{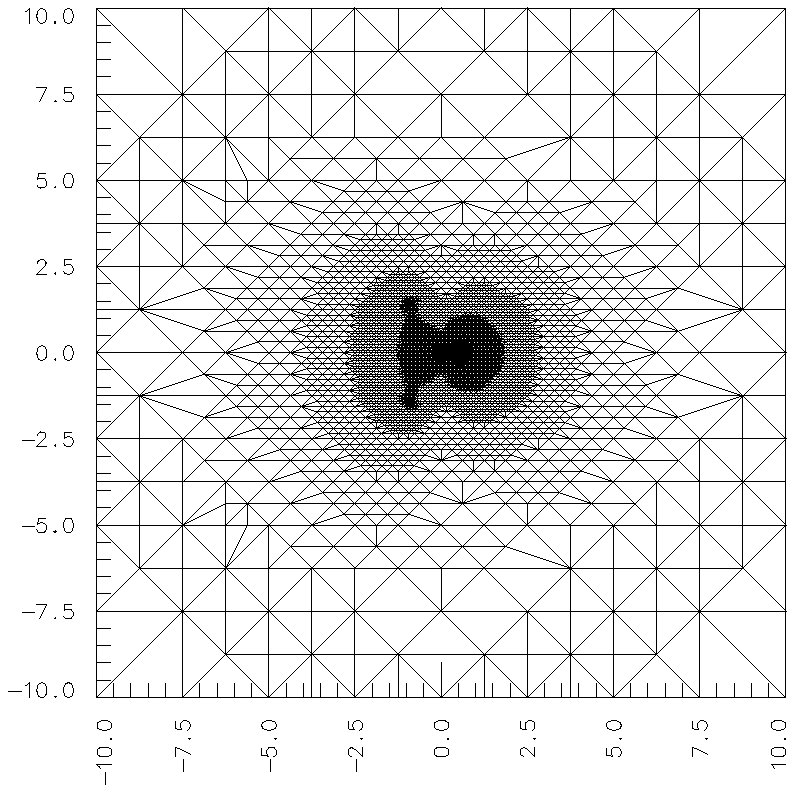}\\
	\includegraphics[width=0.3\linewidth]{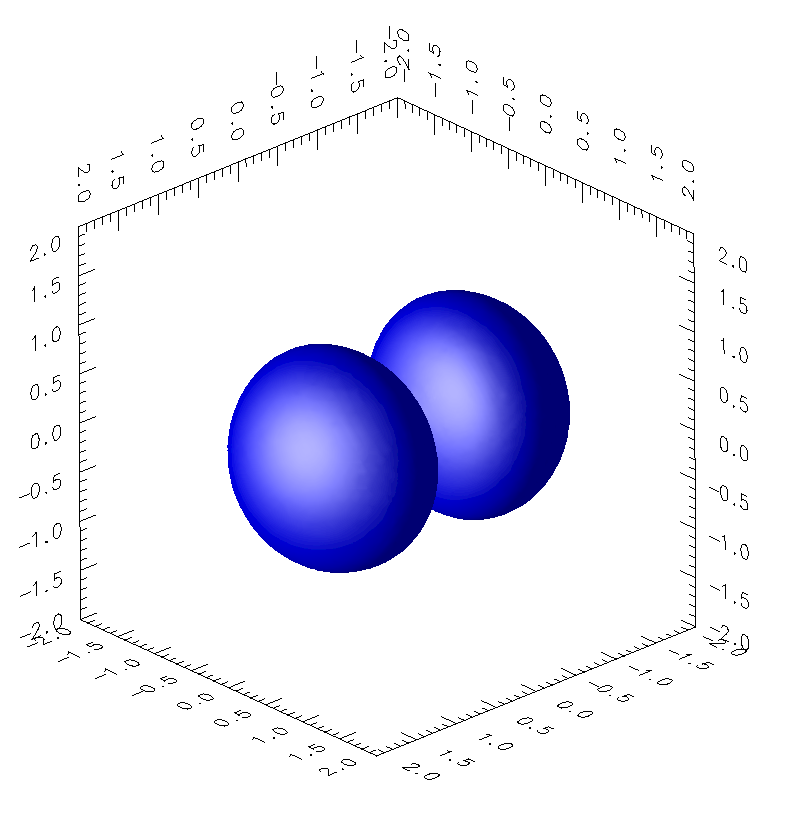}
	\includegraphics[width=0.3\linewidth]{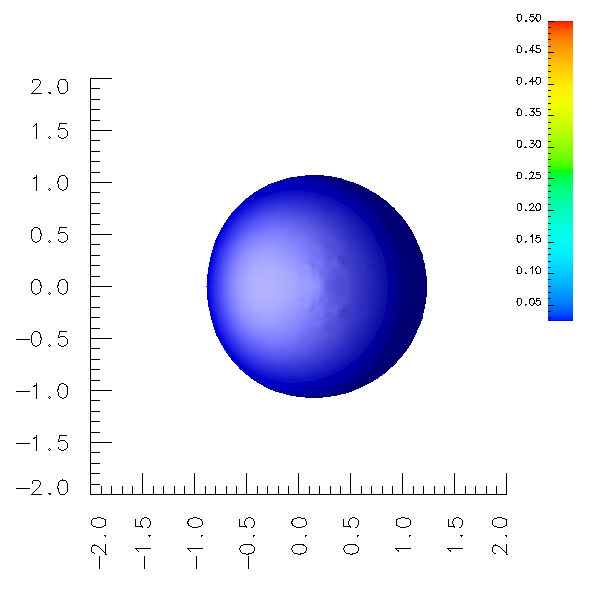}
	\includegraphics[width=0.3\linewidth]{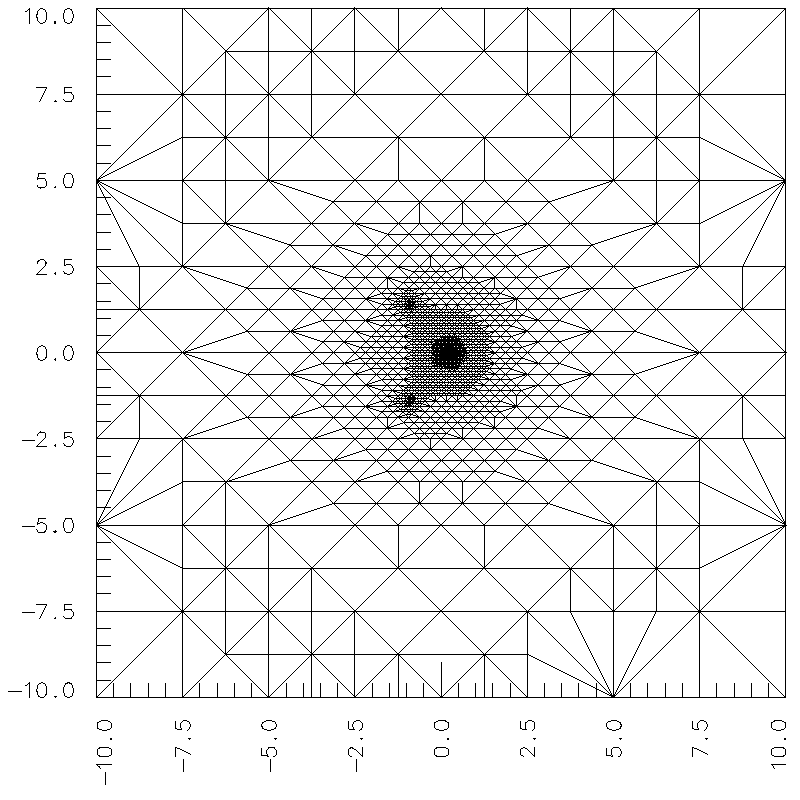}\\
	\includegraphics[width=0.3\linewidth]{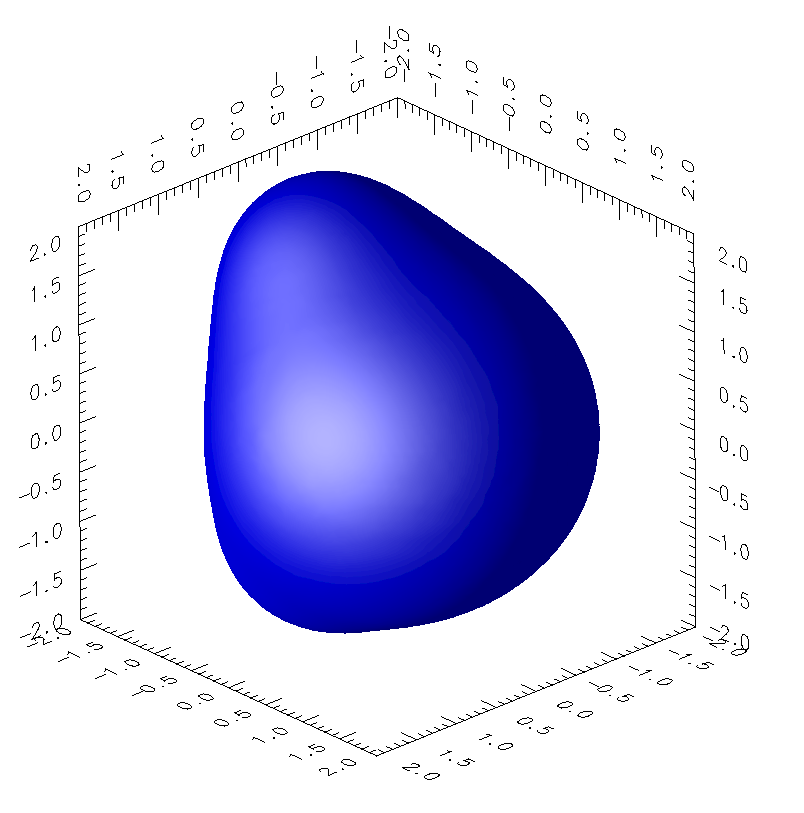}
	\includegraphics[width=0.3\linewidth]{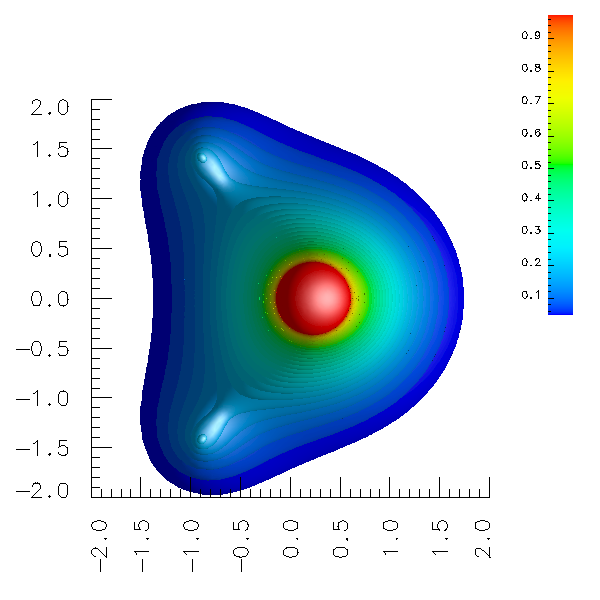}
	\includegraphics[width=0.3\linewidth]{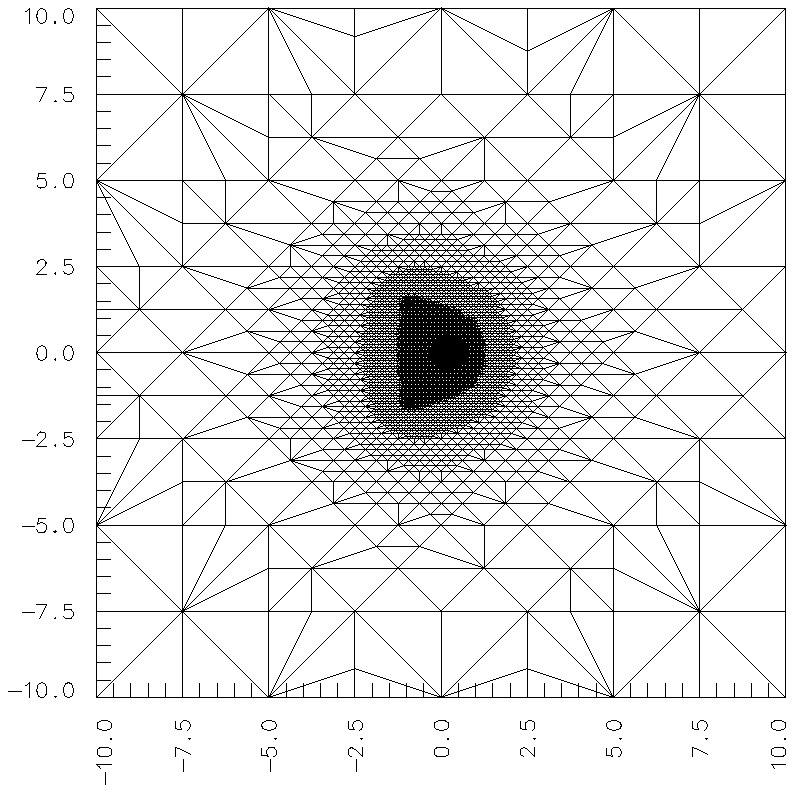}\\
	\includegraphics[width=0.3\linewidth]{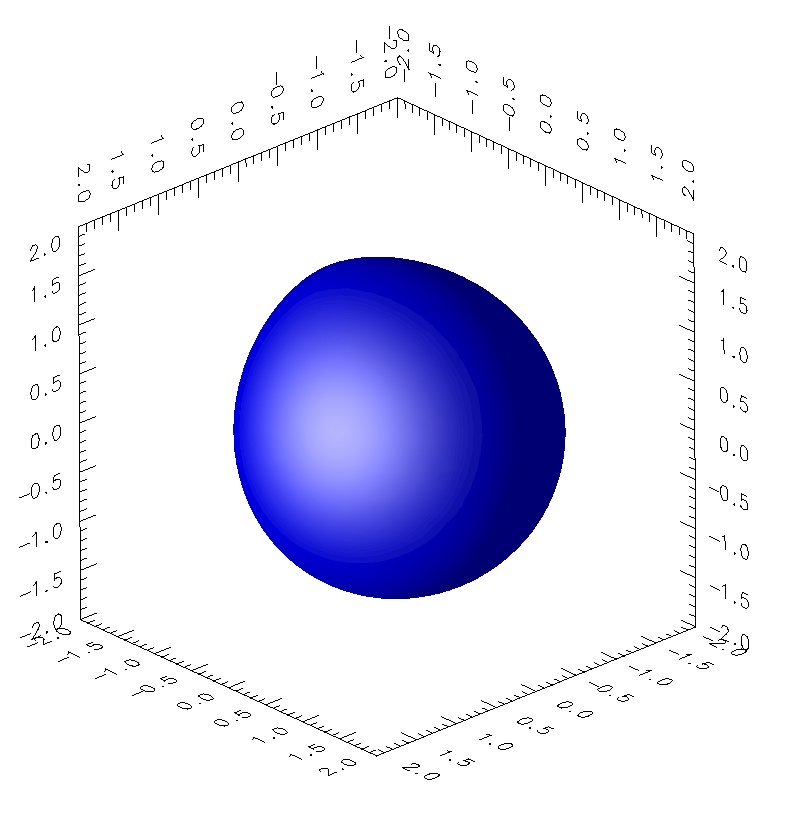}
	\includegraphics[width=0.3\linewidth]{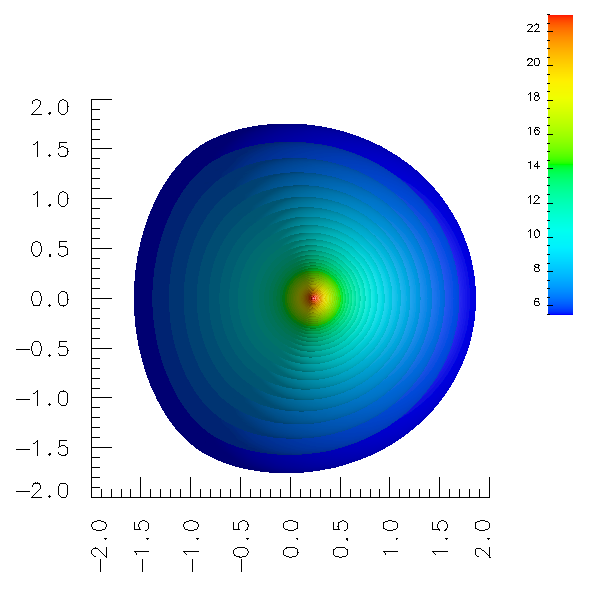}
	\includegraphics[width=0.3\linewidth]{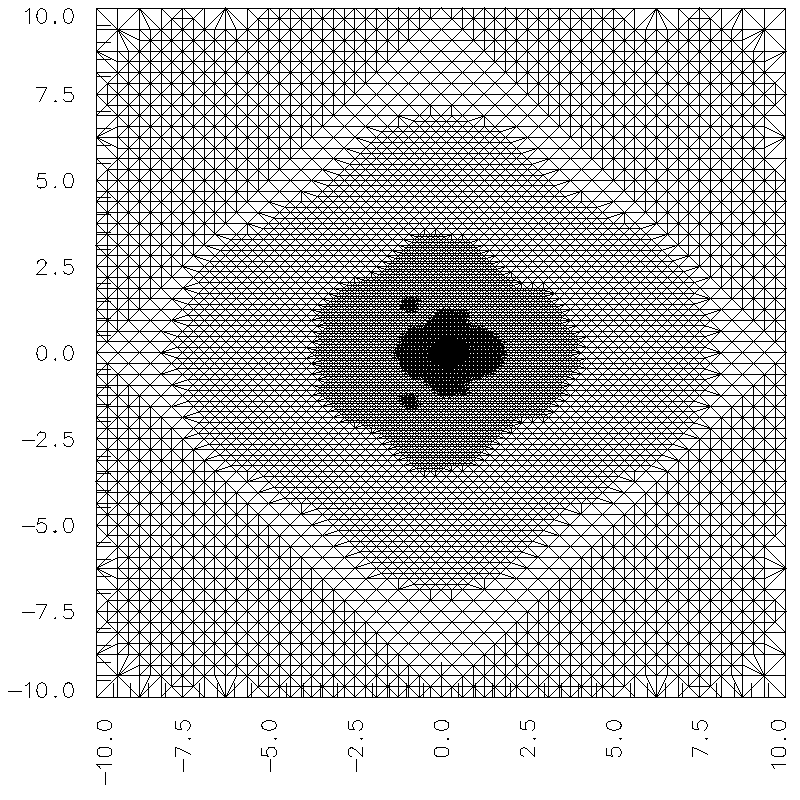}
	\caption{Results for the \ch{H2O} molecule. From top to the fifth row: electron density, sliced electron density on the $X$-$Y$ plane, and mesh distribution for the first, second, third, fourth, and fifth groups of wavefunctions. The sixth row displays the corresponding information for the  electron density on the merged mesh. The last low display the information for the Hartree potential. The number of the mesh grids for these meshes are : 146,653,  155,221, 151,806, 153,808, 152,818, 272,791, 439,968. \label{fig:h2o}}
	
\end{figure}

In the \ch{H2O} example, we first solve its KS equation on a coarse mesh and observe that the 5 eigenvalues are different, then we split the orbitals into 5 group, i.e., one orbital consititute one group. In this way, 7 meshes are generated duiring the whole algorithm, five for the splitted KS orbitals, one is the merged mesh, and the remaining one is for the Hartree potenital. The associated meshes and contours are shown in \Cref{fig:h2o}. The first five rows display the electron density and the mesh distribution for the first, second, third, fourth, and fifth KS orbital, respectively. The sixth row displays the electron density on the merged mesh. The last row displays the information for the Hartree potential. The number of the mesh grids for these meshes are : 146,653, 155,221, 151,806, 153,808, 152,818, 272,791, 439,968.  The results show that the splitting method can effectively separate the wavefunctions into different groups, and the mesh distributions are consistent with the shapes of the electron density and the Hartree potential. The meshes on the right column are consistent with the shape of the electron densities or the Hartree potential. The results shows that the splitting method is feasible in this extreme case, however, the efficiency is not as good as the previous examples. {This can also be read from the splitting factor $sf=0.359$, which is quite small}. The main reason is the difference of the regularities between these five orbitals, especially the last four orbitals is not as evident as the splitting of the core and valence orbitals. Furthermore, the introduction of the new meshes also increases the computational cost, especially for the communication between the meshes.

\subsection{Splitting based on HOMO-LUMO}
It is noted that in density functional theory, the energy $\varepsilon_{\mathrm{HOMO}}$ of the highest occupied KS orbital (HOMO) has physical significance, in the sense that the value $\varepsilon_{\mathrm{HOMO}}$ is in theory equal to the first ionization potential \cite{stowasser1999kohn, fiolhais2008primer}. And the energy $\varepsilon_{\mathrm{LUMO}}$ for the lowest occupied KS orbital (LUMO) is equal to the negative of the electron affinity (EA). Hence, the difference between the HOMO and LUMO energy is the band gap (which is also called HOMO-LUMO gap):
\begin{equation}\label{eq:band-gap}
    E_{\mathrm{gap}} = \varepsilon_{\mathrm{LUMO}}-\varepsilon_{\mathrm{HOMO}}.
\end{equation}
We recognize that the HOMO and LUMO exhibit distinct behaviors, suggesting that a splitting strategy could be effective for this scenario. Specifically, we introduce an additional eigenpair group comprising the LUMO, thereby dividing the problem into three distinct groups: core orbitals, valence orbitals, and the LUMO.

Then this strategy is implemented for the Be atom and present the outcomes in \Cref{tab:be-homo-lumo}. The referenced values are from NIST database \cite{kotochigova2009atomic}. The results indicates that our approach can accurately determine the HOMO-LUMO gap. The three KS meshes utilize 541,117, 540,548, and 766,010 mesh points, respectively, while the merged mesh contains 1,189,514 mesh points. The splitting factor is calculated to be 0.644.

\begin{table*}[!htp]
    \centering
    \begin{tabular}{r|r|r|r|r|r}
    \toprule
    & Reference& split method   & difference & split method (+) & difference \\ \midrule
    $\varepsilon_1$  & -3.856411	&-3.857864	&0.001453	&-3.856486	&0.000075 \\
    $\varepsilon_{\mathrm{HOMO}}$ & -0.205744	&-0.205662	&0.000082	&-0.206129	&0.000385 \\
    $\varepsilon_{\mathrm{LUMO}}$ & -0.077178	&-0.077866	&0.000688	&-0.077331	&0.000153 \\
    $\varepsilon_{\mathrm{gap}}$ & 0.128566	&0.127796	&0.000770	&0.128798	&0.000232
\\
   \bottomrule
    \end{tabular}
    \caption{HOMO-LUMO study for atom Be. Reference values are from NIST database \cite{kotochigova2009atomic}.}
    \label{tab:be-homo-lumo}
\end{table*}

\begin{figure}
	\centering
	\includegraphics[width=0.3\linewidth]{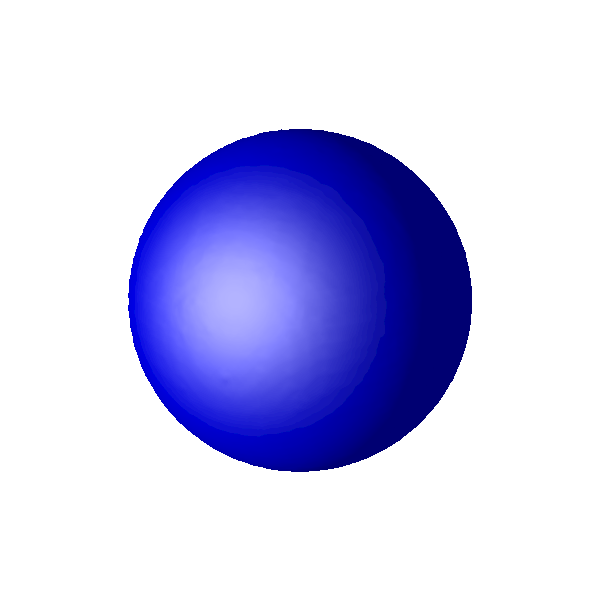}
	\includegraphics[width=0.3\linewidth]{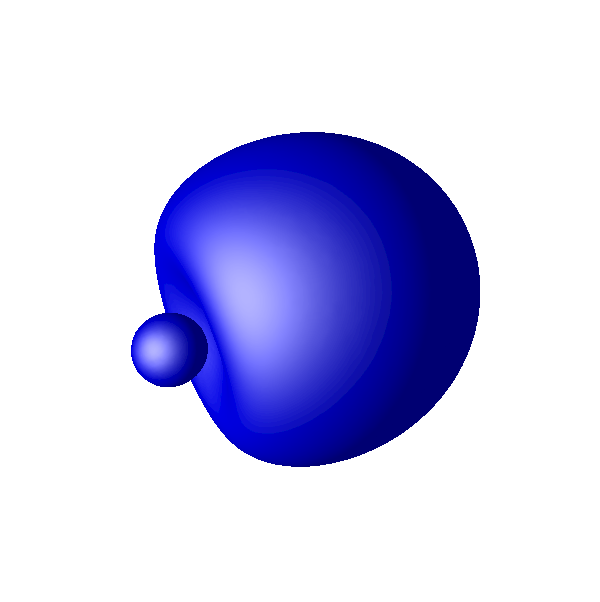}
	\includegraphics[width=0.3\linewidth]{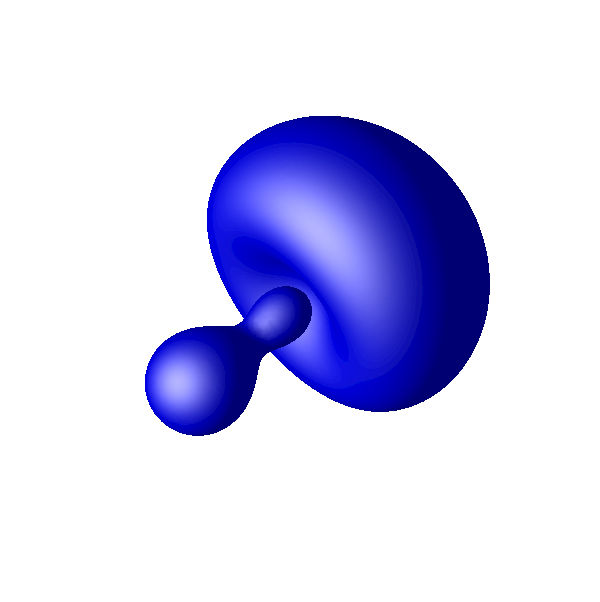}
	\includegraphics[width=0.3\linewidth]{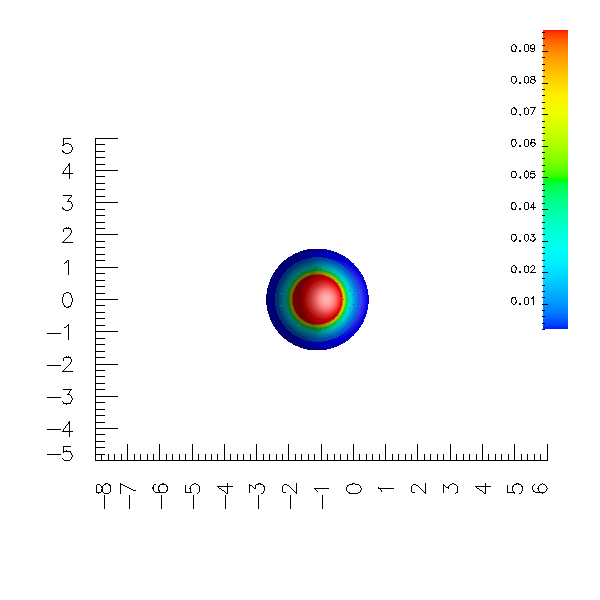}
	\includegraphics[width=0.3\linewidth]{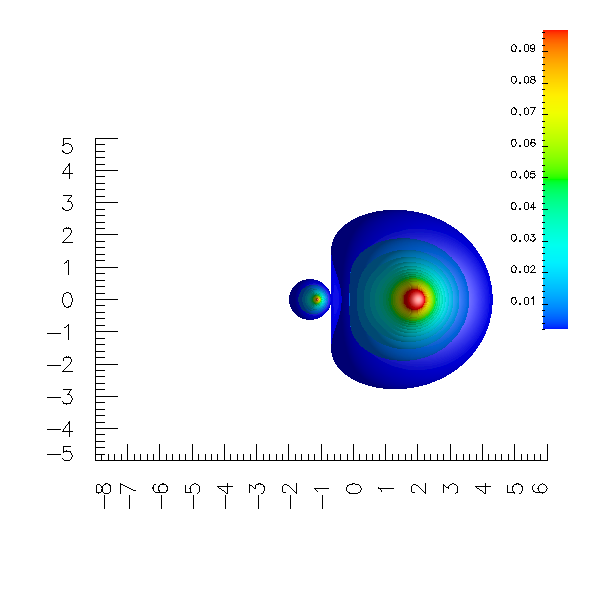}
	\includegraphics[width=0.3\linewidth]{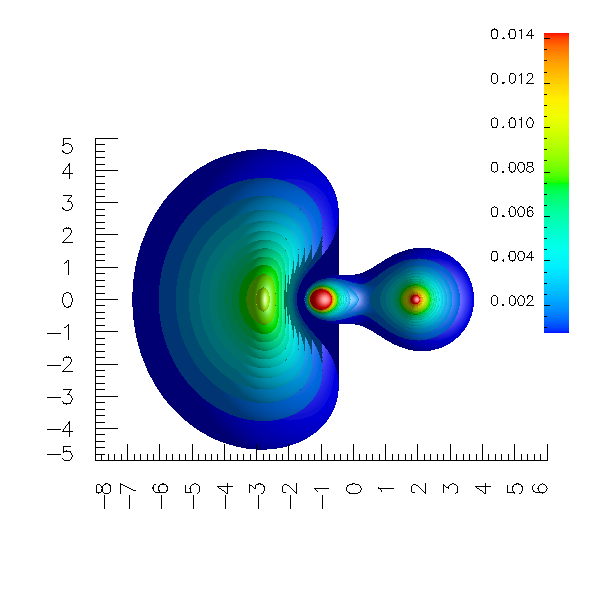}
	\includegraphics[width=0.3\linewidth]{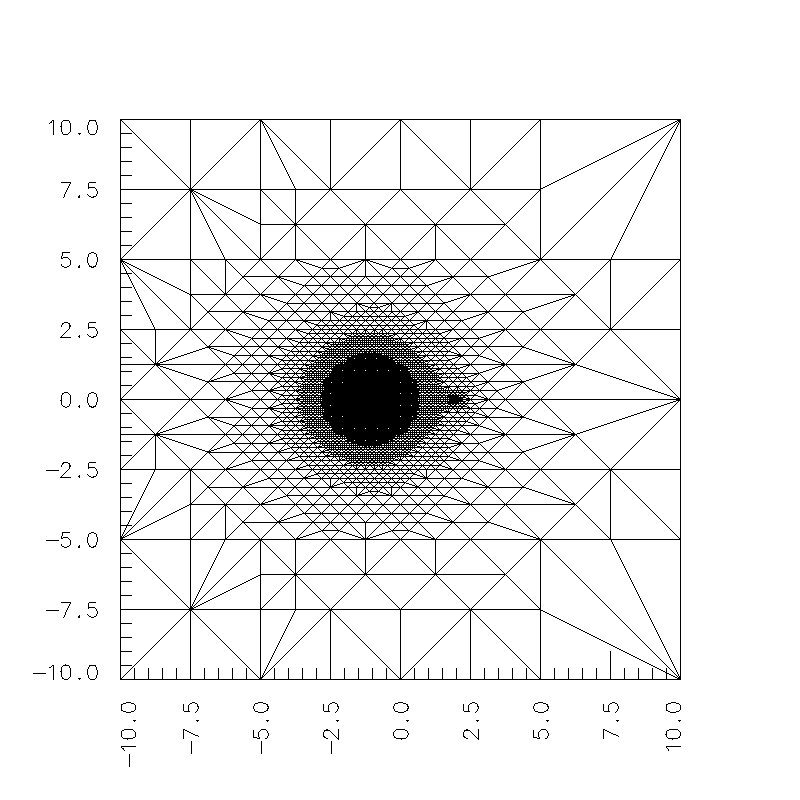}
	\includegraphics[width=0.3\linewidth]{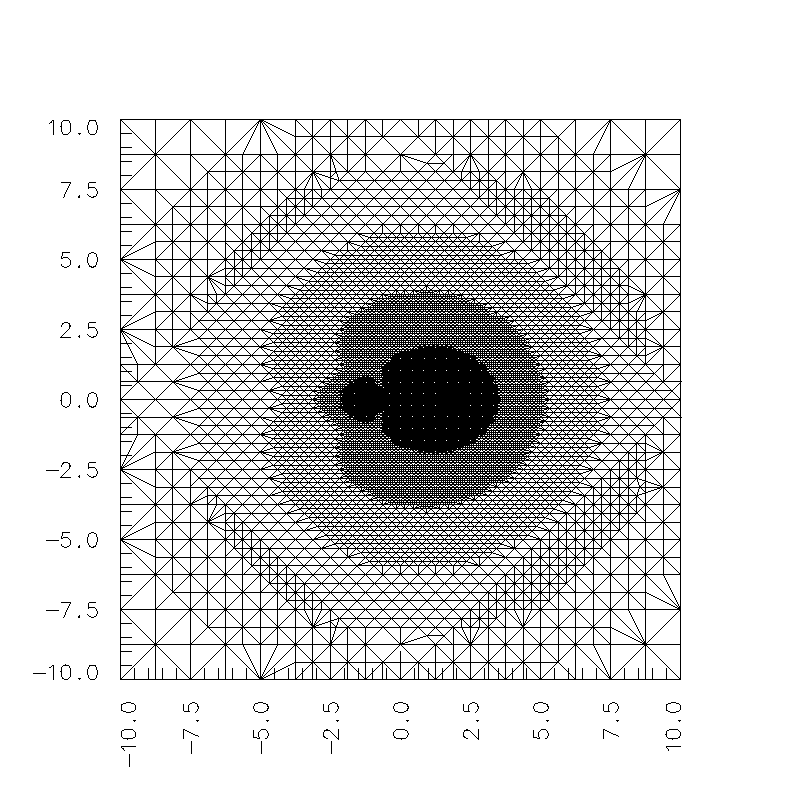}
	\includegraphics[width=0.3\linewidth]{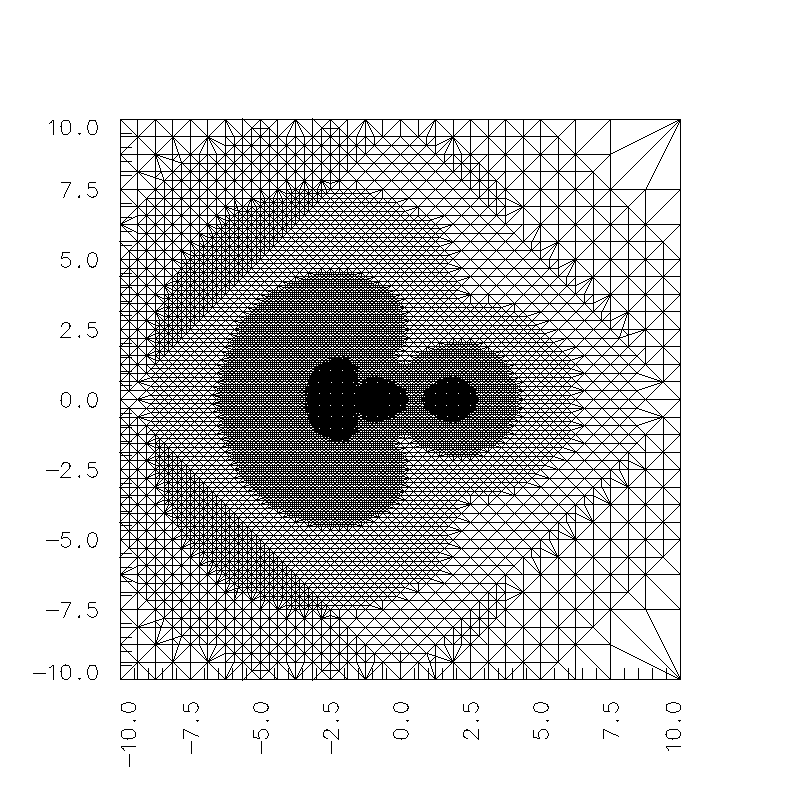}
	\includegraphics[width=0.3\linewidth]{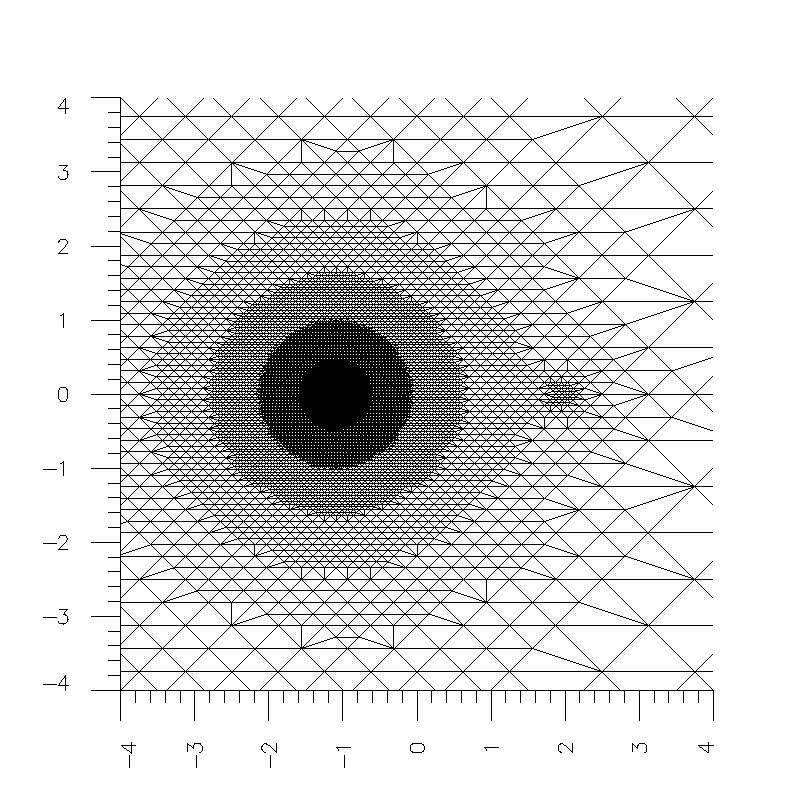}
	\includegraphics[width=0.3\linewidth]{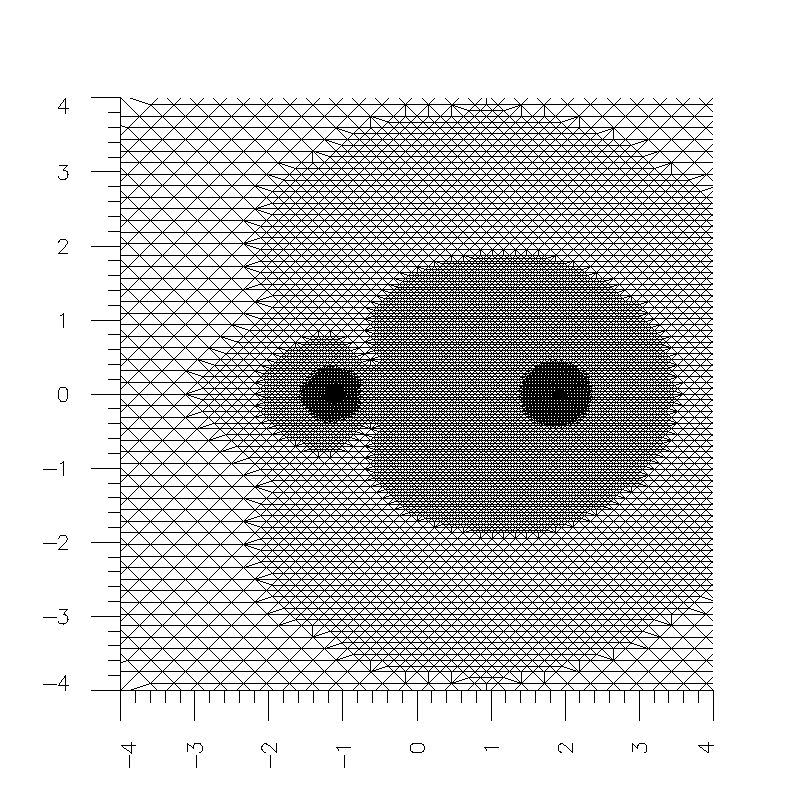}
	\includegraphics[width=0.3\linewidth]{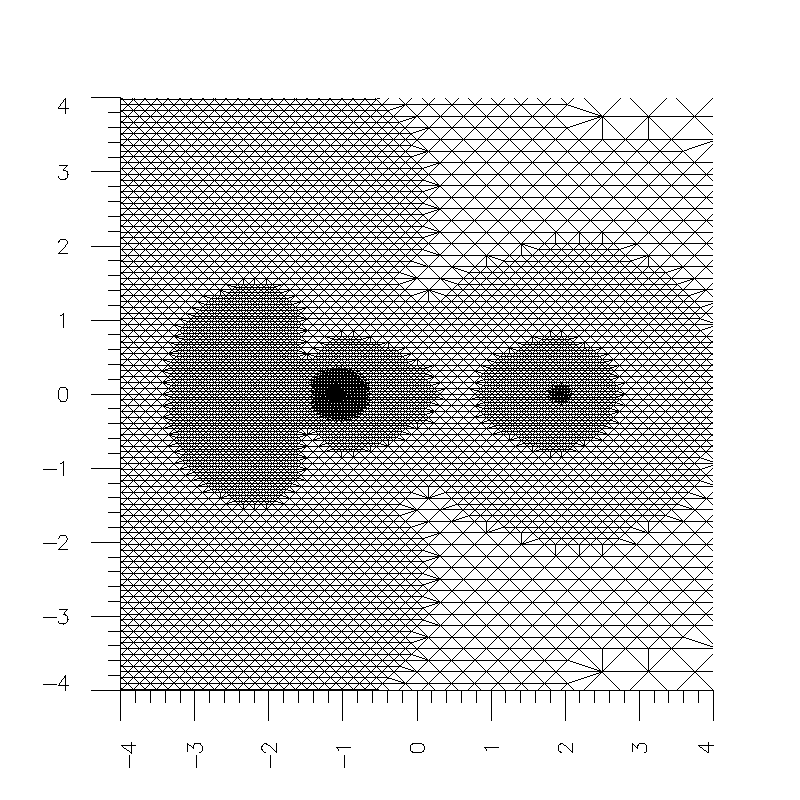}
	\caption{LUMO study for the LiH molecule. \label{fig:lih-lumo}}
\end{figure}

Subsequently, we investigate the HOMO-LUMO for a  LiH molecule example for a clear demonstration, using the splitting strategy based on HOMO-LUMO. The results are depicted in \Cref{fig:lih-lumo}. For the KS meshes, the mesh points are 2,132,630, 1,997,052, and 2,112,120, respectively. The splitting factor is $sf = 0.802$, which illustrate the efficiency of our splitting method.

These experiments demonstrate the feasibility of the HOMO-LUMO-based splitting strategy, offering a potentially efficient approach to deteminating the HOMO-LUMO and calculating the HOMO-LUMO gap.

\section{Conclusion}
In this paper, we presented a novel eigenpair-splitting method for efficiently solving the generalized eigenvalue problem arising from the Kohn-Sham equation. Alternative to traditional DAC methods that rely on domain decomposition, our approach split the problem into several subproblems defined on the whole domain. By integrating a multi-mesh strategy to generate approximate spaces and employing a soft-locking technique to independently solve eigenpairs, the method demonstrates improved computational efficiency. Numerical experiments, including tests on the HOMO-LUMO gap, confirm the effectiveness of the proposed approach. Additionally, we explore strategies for optimal eigenpair grouping. Our discussions indicate that relying solely on eigenvalue-based splitting strategy is not optimal. Instead, we find that strategies focused on core and valence orbitals, as well as those based on the HOMO-LUMO, are more effective.

\section*{Acknowledgement}
The work of Y. Kuang was supported by the National Natural Science Foundation of China (Nos. 12201130, 12326362). The work of G. Hu was supported by the Science and Technology Development Fund, Macao SAR (No. 0068/2024/RIA1), National Natural Science Foundation of China (No. 11922120), MYRG of University of Macau (No. MYRG-GRG2023-00157-FST-UMDF).

\bibliographystyle{unsrt}
\bibliography{Split-KSDFT}

\begin{thebibliography}{10}

\bibitem{kohn1965self}
Walter Kohn and Lu~Jeu Sham.
\newblock Self-consistent equations including exchange and correlation effects.
\newblock {\em Phys. Rev.}, 140(4A):A1133, 1965.

\bibitem{pickett1989pseudopotential}
Warren~E Pickett.
\newblock Pseudopotential methods in condensed matter applications.
\newblock {\em Comput. Phys. Rep.}, 9(3):115--197, 1989.

\bibitem{xiao2010first}
Haiyan~Y Xiao, X~D Jiang, G~Duan, Fei Gao, Xiaotao~T Zu, and William~J Weber.
\newblock {First-principles calculations of pressure-induced phase
  transformation in AlN and GaN}.
\newblock {\em Comput. Mater. Sci.}, 48(4):768--772, 2010.

\bibitem{hamann2013optimized}
DR~Hamann.
\newblock Optimized norm-conserving {V}anderbilt pseudopotentials.
\newblock {\em Phys. Rev. B}, 88(8):085117, 2013.

\bibitem{ainsworth1997posteriori}
Mark Ainsworth and J~Tinsley Oden.
\newblock A posteriori error estimation in finite element analysis.
\newblock {\em Computer methods in applied mechanics and engineering},
  142(1-2):1--88, 1997.

\bibitem{verfurth2013posteriori}
R{\"u}diger Verf{\"u}rth.
\newblock {\em A posteriori error estimation techniques for finite element
  methods}.
\newblock OUP Oxford, 2013.

\bibitem{tsuchida1995electronic}
Eiji Tsuchida and Masaru Tsukada.
\newblock Electronic-structure calculations based on the finite-element method.
\newblock {\em Physical Review B}, 52(8):5573, 1995.

\bibitem{lehtovaara2009all}
Lauri Lehtovaara, Ville Havu, and Martti Puska.
\newblock All-electron density functional theory and time-dependent density
  functional theory with high-order finite elements.
\newblock {\em The Journal of Chemical Physics}, 131(5):054103, 2009.

\bibitem{bao2012h}
Gang Bao, Guanghui Hu, and Di~Liu.
\newblock An $h$-adaptive finite element solver for the calculations of the
  electronic structures.
\newblock {\em Journal of Computational Physics}, 231(14):4967--4979, 2012.

\bibitem{fang2012kohn}
Jun Fang, Xingyu Gao, and Aihui Zhou.
\newblock A {Kohn--Sham} equation solver based on hexahedral finite elements.
\newblock {\em Journal of Computational Physics}, 231(8):3166--3180, 2012.

\bibitem{motamarri2013higher}
Phani Motamarri, Michael~R Nowak, Kenneth Leiter, Jaroslaw Knap, and Vikram
  Gavini.
\newblock Higher-order adaptive finite-element methods for {Kohn--Sham} density
  functional theory.
\newblock {\em Journal of Computational Physics}, 253:308--343, 2013.

\bibitem{chen2014adaptive}
Huajie Chen, Xiaoying Dai, Xingao Gong, Lianhua He, and Aihui Zhou.
\newblock Adaptive finite element approximations for {Kohn--Sham} models.
\newblock {\em Multiscale Modeling \& Simulation}, 12(4):1828--1869, 2014.

\bibitem{maday2014h}
Yvon Maday.
\newblock $h$-{$P$} finite element approximation for full-potential electronic
  structure calculations.
\newblock In {\em Partial Differential Equations: Theory, Control and
  Approximation}, pages 349--377. Springer, 2014.

\bibitem{davydov2016adaptive}
Denis Davydov, Toby~D Young, and Paul Steinmann.
\newblock On the adaptive finite element analysis of the {Kohn--Sham}
  equations: methods, algorithms, and implementation.
\newblock {\em International Journal for Numerical Methods in Engineering},
  106(11):863--888, 2016.

\bibitem{motamarri2020dft}
Phani Motamarri, Sambit Das, Shiva Rudraraju, Krishnendu Ghosh, Denis Davydov,
  and Vikram Gavini.
\newblock {DFT-FE--A massively parallel adaptive finite-element code for
  large-scale density functional theory calculations}.
\newblock {\em Computer Physics Communications}, 246:106853, 2020.

\bibitem{lin2019numerical}
Lin Lin, Jianfeng Lu, and Lexing Ying.
\newblock {Numerical methods for Kohn--Sham density functional theory}.
\newblock {\em Acta Numer.}, 28:405--539, 2019.

\bibitem{dai2020gradient}
Xiaoying Dai, Qiao Wang, and Aihui Zhou.
\newblock {Gradient Flow Based Kohn--Sham Density Functional Theory Model}.
\newblock {\em Multiscale Modeling \& Simulation}, 18(4):1621--1663, 2020.

\bibitem{gao2022orthogonalization}
Bin Gao, Guanghui Hu, Yang Kuang, and Xin Liu.
\newblock An orthogonalization-free parallelizable framework for all-electron
  calculations in density functional theory.
\newblock {\em SIAM Journal on Scientific Computing}, 44(3):B723--B745, 2022.

\bibitem{yang1991direct1}
Weitao Yang.
\newblock Direct calculation of electron density in density-functional theory.
\newblock {\em Physical review letters}, 66(11):1438, 1991.

\bibitem{yang1991direct2}
Weitao Yang.
\newblock Direct calculation of electron density in density-functional theory:
  Implementation for benzene and a tetrapeptide.
\newblock {\em Physical Review A}, 44(11):7823, 1991.

\bibitem{chen2016analysis}
Jingrun Chen and Jianfeng Lu.
\newblock Analysis of the divide-and-conquer method for electronic structure
  calculations.
\newblock {\em Mathematics of Computation}, 85(302):2919--2938, 2016.

\bibitem{li2005multi}
Ruo Li.
\newblock On multi-mesh $h$-adaptive methods.
\newblock {\em Journal of Scientific Computing}, 24(3):321--341, 2005.

\bibitem{kuang2024towards}
Yang Kuang, Yedan Shen, and Guanghui Hu.
\newblock Towards chemical accuracy using a multi-mesh adaptive finite element
  method in all-electron density functional theory.
\newblock {\em Journal of Computational Physics}, 518:113312, 2024.

\bibitem{knyazev2001toward}
Andrew~V Knyazev.
\newblock Toward the optimal preconditioned eigensolver: {L}ocally optimal
  block preconditioned conjugate gradient method.
\newblock {\em SIAM Journal on Scientific Computing}, 23(2):517--541, 2001.

\bibitem{knyazev2007block}
Andrew~V Knyazev, Merico~E Argentati, Ilya Lashuk, and Evgueni~E Ovtchinnikov.
\newblock {Block locally optimal preconditioned eigenvalue xolvers (BLOPEX) in
  HYPRE and PETSc}.
\newblock {\em SIAM Journal on Scientific Computing}, 29(5):2224--2239, 2007.

\bibitem{marques2012libxc}
Miguel A~L Marques, Micael J~T Oliveira, and Tobias Burnus.
\newblock Libxc: {A} library of exchange and correlation functionals for
  density functional theory.
\newblock {\em Computer Physics Communications}, 183(10):2272--2281, 2012.

\bibitem{vosko1980accurate}
Seymour~H Vosko, Leslie Wilk, and Marwan Nusair.
\newblock Accurate spin-dependent electron liquid correlation energies for
  local spin density calculations: a critical analysis.
\newblock {\em Canadian Journal of Physics}, 58(8):1200--1211, 1980.

\bibitem{froese1997computational}
Charlotte Froese-Fischer, Tomas Brage, and Per Jonsson.
\newblock {\em Computational atomic structure: an MCHF approach}.
\newblock Routledge, 1997.

\bibitem{cai2024AFEPack}
Zhenning Cai, Yun Chen, Yana Di, Guanghui Hu, Ruo Li, Wenbin Liu, Heyu Wang,
  Fanyi Yang, Chengbao Yao, and Hongfei Zhan.
\newblock {AFEPack: a general-purpose C++ library for numerical solutions of
  partial differential equations}.
\newblock {\em Communications in Computational Physics}, 2024, in press.

\bibitem{duersch2018robust}
Jed~A Duersch, Meiyue Shao, Chao Yang, and Ming Gu.
\newblock A robust and efficient implementation of {LOBPCG}.
\newblock {\em SIAM Journal on Scientific Computing}, 40(5):C655--C676, 2018.

\bibitem{bao2015real}
Gang Bao, Guanghui Hu, and Di~Liu.
\newblock Real-time adaptive finite element solution of time-dependent
  {Kohn--Sham} equation.
\newblock {\em Journal of Computational Physics}, 281:743--758, 2015.

\bibitem{bao2016towards}
Gang Bao, Guanghui Hu, and Di~Liu.
\newblock Towards translational invariance of total energy with finite element
  methods for {Kohn-Sham} equation.
\newblock {\em Communications in Computational Physics}, 19(1):1--23, 2016.

\bibitem{valiev2010nwchem}
Marat Valiev, Eric~J Bylaska, et~al.
\newblock {NWChem: A comprehensive and scalable open-source solution for large
  scale molecular simulations}.
\newblock {\em Computer Physics Communications}, 181(9):1477--1489, 2010.

\bibitem{ciarlet2002finite}
Philippe~G Ciarlet.
\newblock {\em The finite element method for elliptic problems}.
\newblock SIAM, 2002.

\bibitem{stowasser1999kohn}
Ralf Stowasser and Roald Hoffmann.
\newblock What do the {Kohn-Sham} orbitals and eigenvalues mean?
\newblock {\em Journal of the american chemical society}, 121(14):3414--3420,
  1999.

\bibitem{fiolhais2008primer}
Carlos Fiolhais, Fernando Nogueira, and Miguel~AL Marques.
\newblock {\em A primer in density functional theory}, volume 620.
\newblock Springer, 2008.

\bibitem{kotochigova2009atomic}
Svetlana Kotochigova, Zachary~H Levine, Eric~L Shirley, Mark~D Stiles, and
  Charles~W Clark.
\newblock Atomic reference data for electronic structure calculations (version
  1.6), 2009.
\newblock Available at
  \url{https://www.nist.gov/pml/atomic-reference-data-electronic-structure-calculations/atomic-reference-data-electronic-7},
  Accessed: 2024-08-01.

\end{thebibliography}

\end{document}